\documentclass[aop,MSNbibl,seceqn,nameyear,dvips]{arximspdf}

\doi{10.1214/12-AOP798} 
\volume{41}
\issue{5}
\pubyear{2013}
\firstpage{3658}
\lastpage{3696}

\makeatletter

\newcommand{\blddot}{\bolds{.}}

\newcommand{\cal}{\mathcal}

\newtheorem{theorem}{Theorem}[section]
\newtheorem{lem}{Lemma}[section]
\newtheorem{prop}{Proposition}[section]
\newtheorem{cor}{Corollary}[section]

\newproclaim{Defi}{Definition}[section]
\newproclaim{nota}{Notation}[section]
\newproclaim{rem}{Remark}[section]

\newcommand{\I}{{\mathbf1}}
\newcommand{\bkE}{{\mathbb{E}}}
\newcommand{\p}{{\mathbb{P}}}

\makeatother

\begin{document}
\begin{frontmatter}

\title{Strong approximation results for the empirical process of
stationary sequences}
\runtitle{Strong approximation for the empirical process}

\begin{aug}
\author[A]{\fnms{J\'er\^ome} \snm{Dedecker}\corref{}\ead[label=e1]{jerome.dedecker@parisdescartes.fr}},
\author[B]{\fnms{Florence} \snm{Merlev\`ede}\ead[label=e2]{florence.merlevede@univ-mlv.fr}}
\and
\author[C]{\fnms{Emmanuel} \snm{Rio}\ead[label=e3]{emmanuel.rio@uvsq.fr}}
\runauthor{J. Dedecker, F. Merlev\`ede and E. Rio}
\affiliation{Universit\'e Paris Descartes, Universit\'e Paris-Est
and\break
Universit\'e de Versailles}
\address[A]{J. Dedecker\\
Laboratoire MAP5\\
Universit\'e Paris Descartes\\
Sorbonne Paris Cit\'e\\
UMR 8145 CNRS\\
45 rue des Saints-P\`eres\\
F-75270 Paris cedex 06\\
France\\
\printead{e1}}
\address[B]{F. Merlev\`ede\\
Universit\'e Paris-Est\\
LAMA (UMR 8050)\\
UPEMLV, CNRS, UPEC\\
F-77454 Marne-La-Vall\'ee\\
France\\
\printead{e2}}
\address[C]{E. Rio\\
Laboratoire de math\'ematiques\\
Universit\'e de Versailles\\
UMR 8100 CNRS\\
B\^atiment Fermat\\
45 Avenue des Etats-Unis\\
F-78035 Versailles\\
France\\
\printead{e3}} 
\end{aug}

\received{\smonth{10} \syear{2011}}

%
\begin{abstract}
We prove a strong approximation result for the empirical process
associated to a stationary sequence of real-valued random variables,
under dependence conditions involving only indicators of half lines.
This strong approximation result also holds for the empirical process
associated to iterates of expanding maps with a neutral fixed point at
zero, as soon as the correlations decrease more rapidly than
$n^{-1-\delta}$ for some positive $\delta$. This shows that our
conditions are in some sense optimal.
\end{abstract}

%
\begin{keyword}[class=AMS]
\kwd{60F17}
\kwd{60G10}
\kwd{37E05}
\end{keyword}
\begin{keyword}
\kwd{Strong approximation}
\kwd{Kiefer process}
\kwd{stationary sequences}
\kwd{intermittent maps}
\kwd{weak dependence}
\end{keyword}

\end{frontmatter}

\section{Introduction}
Let $(X_i)_{i \in{\mathbb Z}}$ be a strictly stationary sequence of
real-valued random variables with common distribution function $F$, and
define the empirical process
of $(X_i)_{i \in{\mathbb Z}}$ by
%
\begin{equation}
\label{defempproc} R_{X}(s,t)= \sum_{1 \leq k \leq t}
\bigl({\mathbf1}_{X_k \leq s}-F(s) \bigr),\qquad s \in{\mathbb R},
t \in {\mathbb R}^+.
\end{equation}
For independent identically distributed (i.i.d.) random variables $X_i$
with the uniform distribution over $[0,1]$, \citet{KomMajTus75}
constructed a continuous centered Gaussian process
$K_X$ with covariance function
\[
\bkE\bigl( K_X(s,t)K_X\bigl(s',t'
\bigr) \bigr)= \bigl(t \wedge t'\bigr) \bigl(s\wedge
s' - ss'\bigr)
\]
in such a way that
%
\begin{equation}
\label{strappKMT}\quad \sup_{s \in{\mathbb R}, t \in
[0,1]} \bigl|R_{X}\bigl(s,[nt]
\bigr)-K_X\bigl(s,[nt]\bigr)\bigr| =O\bigl(\log^2 n\bigr)\qquad
\mbox{almost surely}
\end{equation}
[we refer also to \citet{CasLau98} for a detailed
proof]. The rate of convergence given in (\ref{strappKMT}) improves on
the one obtained earlier by \citet{Kie72} and the two-parameter
Gaussian process $K_X$ is known in the literature as the Kiefer process.

Such a strong approximation allows not only to derive weak limit
theorems, as Donsker's invariance principle for the empirical
distribution function, but also almost
sure results, as the functional form of the law of the iterated logarithm
[see \citet{Fin71}]. Moreover, from a statistical point of view,
strong approximations with rates allow to construct many statistical
procedures [we refer to the monograph of \citet{ShoWel86}
which shows how the asymptotic behavior of the empirical process plays
a crucial role in many important statistical applications].

In the dependent setting, the weak limiting behavior of the empirical
process $R_X$ has been studied by many authors in different cases. See,
among many others, the following: \citet{DehTaq89} for stationary
Gaussian sequences, \citet{GirSur02} for linear processes,
\citet{Yu93} for associated sequences,
\citet{BorBurDeh01} for functions of absolutely regular sequences,
\citet{Rio00} for strongly mixing sequences, \citet{Wu08} for
functions of i.i.d. sequences and \citet{Ded10} for
$\beta$-dependent sequences.

Strong approximations of type (\ref{strappKMT}), for the empirical
process with dependent data, have been less studied.
\citet{BerPhi77} proved that, for functions of strongly mixing
sequences satisfying $\alpha(n)=O(n^{-8})$ [where $\alpha(n)$ is the
strong mixing coefficient of \citet{Ros56}], and if $F$ is
continuous, there exists a two-parameter continuous Gaussian process
$K_X$ such that
%
\begin{eqnarray}
\label{strappBP}
&&
\sup_{s \in{\mathbb R}, t \in[0,1]} \bigl|R_{X}\bigl(s,[nt]
\bigr)-K_X\bigl(s,[nt]\bigr)\bigr| \nonumber\\[-8pt]\\[-8pt]
&&\qquad=O\bigl(\sqrt{n} \bigl(\ln(n)
\bigr)^{-\lambda}\bigr) \qquad\mbox{almost surely}\nonumber
\end{eqnarray}
for some $\lambda>0$. The covariance function $\Gamma_X$ of $K_X$ is
given by
\[
\Gamma_X\bigl(s,s',t,t'\bigr) = \min
\bigl(t,t'\bigr)\Lambda_X \bigl(s,s'
\bigr),
\]
where
%
\begin{equation}
\label{covariance} \Lambda_X \bigl(s,s'\bigr)=\sum
_{k \geq0} \operatorname{Cov} ( {\mathbf1}_{X_0
\leq s},
{\mathbf1}_{X_k \leq s'} ) + \sum_{k > 0}
\operatorname{Cov} ( {\mathbf1}_{X_0 \leq s'}, {\mathbf1}_{X_k \leq s} ).
\end{equation}
As a corollary, \citet{BerPhi77} obtained that the sequence
\[
\bigl\{ (2n \ln\ln n)^{-1/2} R_{ X} \bigl(s,[nt]\bigr), n
\geq3 \bigr\}
\]
of random functions on ${\mathbb R} \times[0,1] $ is with probability
one relatively compact for the supremum norm, and
that the set of limit points is the unit ball of the reproducing kernel
Hilbert space (RKHS) associated with $\Gamma_X$. Their result
generalizes the functional form of the Finkelstein's law of the
iterated logarithm. Next, \citet{Yos79} weakened the strong mixing
condition required in \citet{BerPhi77} and proved the strong
approximation (\ref{strappBP}) assuming $\alpha(n)=O(n^{-a})$ for
some $a>3$. However, this condition still appears to be too restrictive:
indeed, Rio [(\citeyear{Rio00}), Theorem 7.2, page 96] proved that the weak
convergence of $n^{-1/2} R_X ( s,n)$ to a Gaussian process holds in
$D({\mathbb R})$ under the weaker condition
$\alpha(n) = O(n^{-a})$ for some $a>1$. In view of this result, one
may think that the strong approximation by a Kiefer process, as given
in (\ref{strappBP}), holds as soon as the dependence coefficients are
of the order of $ O(n^{-a})$ for some $a>1$.

Since the classical mixing coefficients have some limited
applicability, many papers have been written in the last decade to
derive limit theorems under various weak dependence measures [see,
e.g., the monograph by \citet{Dedetal07}]. Concerning the
empirical process, \citet{Ded10} proved that the weak convergence of
$n^{-1/2} R_X ( s,n)$ to a Gaussian process holds in $D({\mathbb R})$
under a dependence condition involving only indicators of a half line,
whereas \citet{Wu08} obtained the same result under conditions on, what he
called, the predictive dependent measures. These predictive dependence
measures allow coupling by independent sequences and are well adapted
to some functions of i.i.d. sequences. However, they seem to be less
adequate for functionals of nonirreducible Markov chains or dynamical
systems having some invariant probability. The recent paper by
\citet{BerHorSch09} deals with strong approximations as in
(\ref{strappBP}) in the weak dependent setting by considering, what
they called, $S$-mixing conditions. Actually, their $S$-mixing
condition lies much closer to the predictive dependent measures
considered by \citet{Wu08} and is also very well adapted to functions of
i.i.d. sequences. Roughly speaking, they obtained (\ref{strappBP}) as
soon as $F$ is Lipschitz continuous, the sequence $(X_i)_{i \in{\mathbb
Z}}$ can be approximated by a $2m$-dependent sequence, and one has a
nice control of the deviation probability of the approximating error.

In this paper, we prove that the strong approximation (\ref{strappBP})
holds under a dependence condition involving only indicators of a half
line, which is quite natural in this context [see the discussion at the
beginning\vspace*{1pt} of Section 2 in \citet{Ded10}]. More
precisely, if $\beta_{2, X}(n)=O(n^{-{(1 + \delta)}})$ for some
positive $\delta$, where the coefficients $\beta_{2, X}(n)$ are defined
in the next section, we prove that there exists a continuous (with
respect to its natural metric) centered Gaussian process $K_X$ with
covariance function given by (\ref{covariance}) such that
%
\begin{equation}
\label{strapp}\quad \sup_{s \in{\mathbb R}, t \in[0,1]} \bigl|R_{X}\bigl(s,[nt]
\bigr)-K_X\bigl(s,[nt]\bigr)\bigr| =O\bigl(n^{1/2 -
\varepsilon}\bigr)\qquad
\mbox{almost surely}
\end{equation}
for some $\varepsilon>0$. As consequences of (\ref{strapp}),
we obtain the functional form of Finkelstein's law of the iterated
logarithm and we recover the empirical central limit theorem obtained
in \citet{Ded10}. Notice that our dependence condition cannot be
directly compared to the one used in the paper by
\citet{BerHorSch09}.

In Theorem~\ref{maindyn} we show that (\ref{strapp}) also holds for
the empirical process associated to an expanding map $T$ of the unit
interval with a neutral fixed point at $0$, as soon as the parameter
$\gamma$
belongs to $]0, 1/2[$ (this parameter describes the behavior of $T$ in
the neighborhood of zero). Moreover, we shall prove that the functional
law of the iterated logarithm
cannot hold
at the boundary $\gamma=1/2$, which shows that our conditions are in
some sense optimal (see Remark~\ref{remoptimality} for a detailed
discussion about the optimality of the conditions).

Let us\vspace*{1pt} now give an outline of the methods used to prove the strong
approximation (\ref{strapp}). We consider the dyadic fluctuations
$ (R_X(s,2^{L+1}) -R_X(s,2^{L}) )_{L \geq0}$ of the empirical
process on a grid with a number of points depending on $L$, let's say
$d_L$. Our proof is mainly based on the existence of multidimensional
Gaussian random variables in ${\mathbb R}^{d_L}$ that approximate, in a
certain sense, the fluctuations of the empirical process on the grid.
These multidimensional Gaussian random variables will be the skeleton
of the approximating Kiefer process. To prove the
existence of these Gaussian random variables, we apply a conditional version
of the Kantorovich--Rubinstein theorem, as given in \citet{Rus85}
(see our Section~\ref{constructkiefer}). The multidimensional Gaussian
random variables are constructed in such a way that the error of
approximation in ${\mathbb L}^1$ of the supremum norm between the
fluctuations of the empirical process on the grid and the
multidimensional Gaussian r.v.'s is exactly the expectation of the
Wasserstein distance of order $1$ (with the distance associated to the
supremum norm) between the conditional law of the fluctuations of the
empirical process on the grid and the corresponding multidimensional
Gaussian law [see Definition~\ref{defwasser} and equality (\ref
{coupling1})]. This error can be evaluated with the help of the
Lindeberg method as done in Section~\ref{sectionGApp} [a~similar
approach has been used recently by \citet{MerRio12} for the
partial sum process]. The oscillations of the empirical process,
namely, the quantities involved in (\ref{psd4}) and (\ref{psd6}), are
handled with the help of a suitable exponential inequality combined
with the Rosenthal-type inequality proved by Dedecker (\citeyear{Ded10}),
Proposition 3.1.
Moreover, it is possible to adapt the method of constructing the
skeleton Kiefer process (by conditioning up to the future rather than
to the past) to deal with the empirical process associated to
intermittent maps.

The paper is organized as follows: in Section~\ref{sectionSS} (resp.,
Section~\ref{sectionIM}) we state the strong approximation results for
the empirical process associated to a class of stationary sequences
(resp., to a class of intermittent maps). Section~\ref{sectionProof} is
devoted to the proof of the main results, whereas some technical tools
are stated and proved in the \hyperref[app]{Appendix}.

\section{Strong approximation for the empirical process associated to
a class of stationary sequences} \label{sectionSS}

Let $(X_i)_{i \in{\mathbb Z}}$ be a strictly stationary sequence of
real-valued random variables defined on the probability space
$(\Omega,{\cal A}, {\mathbb P} )$. Assume that $(\Omega,{\cal A},
{\mathbb P} )$ is large enough to contain a sequence $(U_i)_{i
\in{\mathbb Z}}=(\delta_i, \eta_i)_{i \in{\mathbb Z}}$ of i.i.d. random
variables with uniform distribution over $[0,1]^2$, independent of
$(X_i)_{i \in{\mathbb Z}}$. Define the nondecreasing filtration $({\cal
F}_i)_{i \in{\mathbb Z}}$ by ${\cal F}_i = \sigma(X_k\dvtx  k \leq i)$. Let
${\cal{F}}_{-\infty} = \bigcap_{i \in{\mathbb Z}} {\cal{F}}_{i}$ and
${\cal{F}}_{\infty} = \bigvee_{i \in{\mathbb Z}} {\cal{F}}_{i}$. We
shall denote by ${\mathbb E}_i$ the conditional expectation with
respect to ${\mathcal F}_i$.

Let us now define the dependence coefficients that we consider in this paper.
%
\begin{Defi}\label{beta}
Let $P$ be the law of $X_0$ and $P_{(X_i, X_j)}$ be the law of $(X_i,
X_j)$. Let
$P_{X_k|X_0}$ be the conditional distribution of $X_k$ given $X_0$,
$P_{X_k| {\mathcal
F}_\ell}$ be the conditional distribution of $X_k$ given ${\mathcal
F}_\ell$, and $P_{(X_i,
X_j)|{\mathcal F}_\ell}$ be the conditional distribution of
$(X_i,X_j)$ given ${\mathcal
F}_\ell$. Define the functions $f_t={\mathbf1}_{]-\infty, t]}$, and
$f_t^{(0)}=f_t-P(f_t)$.
Define the random variables
\begin{eqnarray*}
b(X_0, k)&=& \sup_{t \in{\mathbb R}} \bigl|{P}_{X_k|X_0}(f_t)-P(f_t)\bigr|,
\\
b_1({\mathcal F}_{\ell}, k)&=& \sup_{t \in{\mathbb
R}}\bigl|{P}_{X_{k}|{\mathcal F}_\ell}(f_t)-P(f_t)\bigr|,
\\
b_2({\mathcal F}_{\ell}, i,j)&=& \sup_{(s,t) \in{\mathbb
R}^2}\bigl|P_{(X_i, X_j)|{\mathcal F}_\ell}
\bigl(f_t^{(0)}\otimes f_s^{(0)}
\bigr)-P_{(X_i, X_j)}\bigl(f_t^{(0)}\otimes
f_s^{(0)}\bigr)\bigr|.
\end{eqnarray*}
Define now the coefficients
\[
\beta\bigl(\sigma(X_0), X_k\bigr)={\mathbb E}
\bigl(b(X_0, k)\bigr),\qquad \beta_{1,X}(k)= {\mathbb E}
\bigl(b_1({\mathcal F}_0, k)\bigr)
\]
and
\[
\beta_{2,X}(k) =\max\Bigl\{\beta_1(k),
\sup_{i>j\geq k} {\mathbb E}\bigl(\bigl(b_2({\mathcal
F}_0, i,j)\bigr)\bigr)\Bigr\}.
\]
Define also
\[
\alpha_{1,X}(k)=\sup_{t \in{\mathbb R}} \bigl\|{P}_{X_{k}|{\mathcal
F}_0}(f_t)-P(f_t)
\bigr\|_1
\]
and note that $\alpha_{1, X}(k) \leq\beta_{1, X}(k) \leq\beta_{2, X}(k)$.
\end{Defi}

Examples of nonmixing sequences $(X_i)_{i \in{\mathbb Z}}$ in the
sense of \citet{Ros56} for which the coefficients $\beta_{2,X}(n)$
can be computed
may be found in the
paper by \citet{DedPri07}. Let us give a first elementary
example. Let $X_i= \sum_{k\geq0}
a_k \varepsilon_{i-k}$, where $(\varepsilon_i)_{i \in{\mathbb Z}}$
is a sequence of
i.i.d. random variables such that ${\mathbb E}(|\varepsilon_0|^\alpha) <
\infty$ for some
$\alpha>0$, and $a_i=O(\rho^i)$ for some $\rho\in{]0, 1[}$. Let $w$
be the
modulus of continuity of $F$. If
\[
w(x)\leq C \bigl|\ln(x)\bigr|^{-a}\qquad \mbox{in a neighborhood of 0, for some
$a>1$},
\]
then $\beta_{2, X}(n) =O(n^{-a})$ [see Remark 2.3 in \citet{Ded10}].
We shall present another example in the next section.

Our main result is the following:
%
\begin{theorem}\label{mainemp}
Assume that $\beta_{2,X}(n)=O(n^{-1-\delta})$ for some $\delta>0$.
Then:\vadjust{\goodbreak}

\begin{longlist}[(2)]
\item[(1)] For all
$(s,s')$ in ${\mathbb R}^2$, the series $\Lambda_X (s,s')$ defined by
(\ref{covariance}) converges absolutely.
%
\item[(2)] For any $(s,s')\in{\mathbb R}^2$ and $(t,t')$ in ${\mathbb R}^+
\times{\mathbb R}^+$, let $\Gamma_X(s,s',t,t') =\break \min(t,t')\*\Lambda_X
(s,s')$. There exists a centered Gaussian process $K_X$ with
covariance function~$\Gamma_X$, whose sample paths
are almost surely uniformly continuous with respect to the pseudometric
\[
d\bigl((s,t),\bigl(s',t'\bigr)\bigr)=\bigl|F(s)-F
\bigl(s'\bigr)\bigr|+\bigl|t-t'\bigr|
\]
and such that (\ref{strapp}) holds with
$\varepsilon= \delta^2 / (22 ( \delta+2)^2 )$.
\end{longlist}
\end{theorem}

Note that we do not make any assumption on the continuity of the distribution
function $F$.

As in the paper of \citet{BerHorSch09}, we can
formulate corollaries to Theorem~\ref{mainemp}. The first one is direct.
Let $D(\mathbb R \times[0,1])$ be the Skorohod space equipped with the
Skorohod topology, as described in \citet{BicWic71}.
%
\begin{cor} \label{cor1} Assume that $\beta_{2,X}(n)=O(n^{-1-\delta
})$ for some $\delta>0$. Then the empirical process $\{ n^{-1/2} R_{
X} (s,[nt]), s\in{\mathbb R}, t \in[0,1]\}$ converges in $D(\mathbb
R \times[0,1])$ to the Gaussian process $K_X$
defined in item (2) of Theorem~\ref{mainemp}.
\end{cor}

To obtain the second one, we need to combine the strong approximation
(\ref{strapp}) with Theorem 2 in \citet{Lai74}.
%
\begin{cor} \label{cor2} Assume that $\beta_{2,X}(n)=O(n^{-1-\delta
})$ for some $\delta>0$. Then, with probability one, the sequence $\{
(2n \ln\ln n)^{-1/2} R_{ X} (s,[nt]), n \geq3 \}$ of random functions
on ${\mathbb R} \times[0,1] $ is relatively compact for the supremum
norm, and the set of limit points is the unit ball of the reproducing
kernel Hilbert space (RKHS) associated with the covariance
function $\Gamma_X$ defined in Theorem~\ref{mainemp}.
\end{cor}

\section{Strong approximation for the empirical process associated to
a class of intermittent maps} \label{sectionIM}

In this section we consider the following class of intermittent maps,
introduced in \citet{DedGouMer10}:
%
\begin{Defi}
A map $T\dvtx [0,1] \to[0,1]$ is a generalized Pomeau--Man\-neville
map (or GPM map) of parameter $\gamma\in{]0,1[}$ if there exist
$0=y_0<y_1<\cdots<y_d=1$ such that, writing $I_k={]y_k,y_{k+1}[}$,
\begin{longlist}[(4)]
\item[(1)] The restriction of $T$ to $I_k$ admits a $C^1$ extension
$T_{(k)}$ to $\overline{I}_k$.
\item[(2)] For $k\geq1$, $T_{(k)}$ is $C^2$ on $\overline{I}_k$, and $\inf
_{x \in\overline{I}_k}|T_{(k)}'(x)|>1$.
\item[(3)] $T_{(0)}$ is $C^2$ on $]0, y_1]$, with $T_{(0)}'(x)>1$ for $x\in
(0,y_1]$, $T_{(0)}'(0)=1$ and $T_{(0)}''(x) \sim c
x^{\gamma-1}$ when $x\to0$, for some $c>0$.
\item[(4)] $T$ is topologically transitive, that is, there exists some $x$
in $]0,1[$ such that $\{ T^n (x)\dvtx  n \in{\mathbb N} \}$ is a
dense subset of $]0,1[$.\vadjust{\goodbreak}
\end{longlist}
\end{Defi}
The third condition ensures that $0$ is a neutral fixed point
of $T$, with $T(x)=x+c' x^{1+\gamma} (1+o(1))$ when $x\to0$.
The fourth condition is necessary to avoid situations where
there are several absolutely continuous invariant measures or
where the neutral fixed point does not belong to the support of
the absolutely continuous invariant measure.
As a well-known example of a GPM map, let us cite the
\citet{LivSauVai99} map (LSV map) defined by
\[
T(x)= %
\cases{ x\bigl(1+ 2^\gamma x^\gamma
\bigr), &\quad if $x \in[0, 1/2]$,
\cr
2x-1, &\quad if $x \in(1/2, 1]$.}
\]
Theorem 1 in \citet{Zwe98} shows that a GPM map $T$ admits a
unique absolutely continuous invariant probability measure
$\nu$, with density $h_\nu$. Moreover, it is ergodic, has full
support, and $h_\nu(x)/x^{-\gamma}$ is bounded from above and
below.

Let $Q$ be the Perron--Frobenius operator of $T$ with respect to $\nu$,
defined~by
%
\begin{equation}
\label{Perron} \nu(f \cdot g\circ T)=\nu\bigl(Q(f) g\bigr)
\end{equation}
for any bounded measurable functions $f$ and $ g$. Let $(X_i)_{i \in
{\mathbb Z}}$ be a stationary Markov chain with invariant measure $\nu$
and transition Kernel $Q$. Dedecker and Prieur [(\citeyear{DedPri09}),
Theorem 3.1] have proved that
%
\begin{equation}
\label{boundonbeta} \beta_{2, X}(n)= O\bigl(n^{-a}\bigr)\qquad
\mbox{for any $a< (1 - \gamma)/\gamma$}
\end{equation}
[this upper bound
was stated for the Liverani--Saussol--Vaienti
map only, but is also valid in our context: see the last paragraph of
the introduction
in \citet{DedPri09}]. As a consequence, if $\gamma< 1/2$,
the stationary sequence
$(X_i)_{i \in{\mathbb Z}}$ satisfies all the assumptions of Theorem
\ref{mainemp}.

Now $(T, T^2,\ldots, T^n)$ is distributed as
$(X_n,X_{n-1},\ldots, X_1)$ on $([0, 1], \nu)$ [see, e.g.,
Lemma XI.3 in \citet{HenHer01}].
Hence, any information on the law of the sums
$\sum_{i=1}^n (f \circ T^i-\nu(f))$ can be
obtained by studying the law of $\sum_{i=1}^n (f(X_i)-\nu(f))$.
However, the reverse time property cannot be used directly to transfer
the almost sure results for $\sum_{i=1}^n (f(X_i)-\nu(f))$ to the sum
$\sum_{i=1}^n (f \circ T^i-\nu(f))$.

For any $s \in[0,1]$ and $t\in{\mathbb R}$, let us consider the
empirical process associated to the dynamical system $T$:
%
\begin{equation}
\label{defempdyn} R_{T}(s,t)= \sum_{1 \leq i \leq t}
\bigl( {\mathbf1}_{T^i \leq s} -F_{\nu} (s)\bigr) \qquad\mbox{where
}F_{\nu} (s) = \nu\bigl([0,s]\bigr).
\end{equation}

For any $\nu$-integrable function $g$, let $g^{(0)}=g-\nu(g)$ and
recall that $f_s={\mathbf1}_{]-\infty, s]}$.
Our main result is the following:
%
\begin{theorem}\label{maindyn}
Let $T$ be a GPM map with parameter $\gamma\in{]0, 1/2[}$. Then:
\begin{longlist}[(2)]
\item[(1)]
For all $(s,s')\in[0,1]^2$, the following series
converges absolutely:
%
\begin{equation}
\label{defcovdyn} \Lambda_{T} \bigl(s,s'\bigr)= \sum
_{k \geq0} \nu\bigl(f_s^{(0)}\cdot
f_{s'}^{(0)}\circ T^k\bigr) + \sum
_{k >0} \nu\bigl(f_{s'}^{(0)}\cdot
f_{s}^{(0)}\circ T^k\bigr).\vadjust{\goodbreak}
\end{equation}
\item[(2)]
For any $(s,s')\in[0,1]^2$ and any $(t,t')\in{\mathbb R}^+ \times
{\mathbb R}^+$, let $\Gamma_{T}(s,s',t,t') = \min(t,t')\Lambda_{T}
(s,s')$. There exists a continuous centered Gaussian process $K^*_T$
with covariance function $\Gamma_T$ such that for some $\varepsilon>0$,
\[
\sup_{(s,t) \in[0,1]^2} \bigl|R_{T}\bigl(s,[nt]\bigr)-K^*_T
\bigl(s,[nt]\bigr)\bigr| =O\bigl(n^{1/2 -
\varepsilon}\bigr) \qquad\mbox{almost surely}.
\]
\end{longlist}
\end{theorem}
%
\begin{rem}
According to the proof of Theorem~\ref{maindyn}, item (2) holds for any
$\varepsilon$ in
$]0, (1-2 \gamma)^2/ 22 [$.
\end{rem}
%
\begin{rem} \label{remoptimality} In the case $\gamma=1/2$, Dedecker
[(\citeyear{Ded10}), Proposition 4.1] proved that, for the LSV map with $\gamma
=1/2$, the finite-dimensional marginals of the process $\{(n \ln n
)^{-1/2}R_T(\cdot, n)\}$ converge in distribution to those of the
degenerated Gaussian
process $G$ defined by
\[
\mbox{for any $t \in[0, 1]$}\qquad G(t)= \sqrt{h_\nu(1/2)}
\bigl(1-F_\nu(t)\bigr){\mathbf{1}}_{t \neq0} Z,
\]
where $Z$ is a standard normal.
This shows that an approximation by a Kiefer process as in Theorem \ref
{maindyn} cannot hold at the boundary $\gamma=1/2$.

For the same reason, when $\gamma=1/2$, the conclusion of Theorem \ref
{mainemp} does not apply to the stationary Markov chain $(X_i)_{i \in
{\mathbb Z}}$ with invariant measure $\nu$ and transition kernel $Q$
given in (\ref{Perron}).
In fact, it follows from Theorem 3.1 in \citet{DedPri09} that
$\beta_{2, X}(k)>C/k$ for some positive constant $C$, so that the
Markov chain $(X_i)_{i \in{\mathbb Z}}$ does not satisfy the assumptions
of Theorem~\ref{mainemp}.

In the case $\gamma= 1/2$, with the same proof as that of Theorem 1.7
of \citet{DedGouMer10}, we see
that, for any $(s,t) \in[0, 1]^2$ and $b>1/2$,
\[
\lim_{n \rightarrow\infty} \frac{1}{\sqrt n (\ln n )^b} R_{T}\bigl
(s,[nt]\bigr)=0
\qquad\mbox{almost everywhere.}
\]
This almost sure result is of the same flavor as in the corresponding
i.i.d. case, when the random variables have exactly a weak moment of order
2, so that the normalization in the central
limit theorem is $(n \ln n )^{-1/2}$: see the discussion in
\citet{DedGouMer10}, last paragraph of Section 1.2.
\end{rem}

\section{Proofs} \label{sectionProof}

In this section we shall sometimes use the notation $a_{n}\ll b_{n}$ to
mean that there exists a numerical constant $C$
not depending on $n$ such that $a_{n}\leq C b_{n}$, for all positive
integers $n$.

\subsection{\texorpdfstring{Proof of Theorem \protect\ref{mainemp}}{Proof of Theorem 2.1}}\label{sec4.1}
Notice first that for any $(s,s') \in{\mathbb R}^2$,
\[
\bigl| \operatorname{Cov} ( \I_{X_0 \leq s},\I_{X_k \leq s'} )\bigr| \leq\bigl\|
\bkE_0 \bigl( \I_{X_k \leq s'} - F\bigl(s'\bigr)
\bigr) \I_{X_0 \leq s} \bigr\|_1 \leq\bkE\bigl(b(X_0,k)
\bigr) \leq\beta_{1,X} (k).
\]
Since $\sum_{k \geq0} \beta_{1,X} (k) < \infty$, item (1) of Theorem
\ref{mainemp} follows.

To prove item (2), we first introduce another probability on $\Omega$.
Let ${\mathbb P}_0^*$ be the probability on $\Omega$ whose density
with respect to
${\mathbb P}$ is
%
\begin{equation}
\label{defnewProb}\qquad C(\beta)^{-1} \Biggl( 1+4 \sum _{k=1}^\infty b(X_0,
k) \Biggr) \qquad\mbox{with }C(\beta) = 1+4 \sum_{k=1}^\infty\beta
\bigl(\sigma(X_0), X_k\bigr).
\end{equation}
Recall that $P$ is the distribution of $X_0$. Then the image measure $P^*$
of ${\mathbb P}_0^*$ by $X_0$
is absolutely continuous with respect to $P$
with
density
%
\begin{equation}
\label{defnewprob} C(\beta)^{-1} \Biggl(1+4 \sum
_{k=1}^\infty b(x, k) \Biggr).
\end{equation}
Let $F_{P^*}$ be the distribution function of $P^*$, and
let $F_{P^*}(x-0)= \sup_{z<x} F_{P^*}(z)$.
Recall that the sequence $(\eta_i)_{i \in{\mathbb Z}}$ of i.i.d. random
variables with uniform distribution over $[0,1]$ has been introduced
at the beginning of Section~\ref{sectionSS}.
Define then the random variables
%
\begin{equation}
\label{newrv} Y_i=F_{P^*}(X_i-0)+
\eta_i \bigl(F_{P^*}(X_i)-F_{P^*}(X_i-0)
\bigr).
\end{equation}

Let $P_Y$ be the distribution of $Y_0$ and $F_Y$ be the distribution
function of $Y_0$.
Some properties of the sequence $(Y_i)_{i \in{\mathbb Z}}$ are
given in Lemma~\ref{propY} of the \hyperref[app]{Appendix}. In particular, it follows
from Lemma
\ref{propY}
that $X_i=F_{P^*}^{-1}(Y_i)$ almost surely, where $F_{P^*}^{-1}$ is the
generalized inverse of the cadlag function $F_{P^*}$.
Hence,
$R_{ X}(\cdot, \cdot)= R_{ Y} (F_{P^*}(\cdot), \cdot)$ almost
surely, where
\[
R_{ Y}(s,t)= \sum_{1 \leq k \leq t} \bigl(
\I_{Y_k \leq s}-F_Y(s) \bigr),\qquad s \in{[0, 1]}, t \in{\mathbb
R}^+.
\]
%

We now prove that, if $\beta_{2, X}(n)=O(n^{-1-\delta})$ for some
$\delta>0$, then
the conclusion of
Theorem~\ref{mainemp} holds for the stationary sequence $(Y_i)_{i \in
{\mathbb Z}}$ and the associated continuous Gaussian process $K_Y$ with
covariance function $\Gamma_Y(s,s',t,t') = \min(t,t')\Lambda_Y
(s,s')$, where
%
\begin{equation}
\label{defgammaY} \Lambda_Y \bigl(s,s'\bigr)=\sum
_{k \geq0} \operatorname{Cov} ( \I_{Y_0 \leq
s},
\I_{Y_k \leq s'} ) + \sum_{k > 0} \operatorname{Cov}
( \I_{Y_0
\leq s'},\I_{Y_k \leq s} ).
\end{equation}
This implies Theorem~\ref{mainemp}, since
$\Gamma_X(s,s', t, t')= \Gamma_Y(F_{P^*}(s), F_{P^*}(s'), t, t')$.

The proof is divided in two steps: the construction of the Kiefer
process with the help of a conditional version of the
Kantorovich--Rubinstein theorem and a probabilistic upper bound for the
error of approximation.

\subsubsection{Construction of the Kiefer process} \label{constructkiefer}

For $L \in{\mathbb N}$, let $m(L) \in{\mathbb N}$ and $r(L) \in
{\mathbb N}^*$ be such that $ m (L) \leq L$ and $4r(L)\leq m(L)$. For
$j $ in $ \{1,\ldots, 2^{r(L)} -1\} $, let $s_j = j 2^{-r(L)}$ and
define for any $\ell\in\{1,\ldots, 2^{L-m(L)} \} $,
\[
I_{L, \ell} = {\bigl]2^L + (\ell-1)2^{m(L)},
2^L + \ell2^{m(L)}\bigr]} \cap {\mathbb N}
\]
and
\[
U^{(j)}_{L,\ell} = \sum_{i\in
I_{L, \ell}} \bigl(
\I_{Y_i \leq s_j}-F_Y(s_j) \bigr).
\]
The associated column vectors $U_{L, \ell}$ are then defined in
${\mathbb R}^{2^{r(L)} -1}$ by
\[
U_{L, \ell} = \bigl(U^{(1)}_{L,\ell},\ldots,
U^{(2^{r(L)}
-1)}_{L,\ell} \bigr)'.
\]
Let us now introduce some definitions.
%
\begin{Defi} \label{defwasser} Let $m$ be a positive integer. Let
$P_1$ and $P_2$ be two probabilities on $({\mathbb R}^{m}, {\mathcal
B}({\mathbb R}^{m}))$. Let $d$ be a distance on ${\mathbb R}^{m}$
associated to a norm. The Wasserstein distance of order $1$ between
$P_1$ and $P_2$ with respect to the distance $d$ is defined by
\begin{eqnarray*}
W_{d} (P_1, P_2) &=& \inf\bigl\{ {\mathbb E}
\bigl(d(X,Y)\bigr), (X,Y) \mbox{ such that }X \sim P_1, Y \sim
P_2 \bigr\}\\
&=& \sup_{f \in\operatorname{Lip} (d) } \bigl( P_1(f) -
P_2(f)\bigr),
\end{eqnarray*}
where $\operatorname{Lip}(d) $ is the set of functions from ${\mathbb R}^m$
into ${\mathbb R}$ that are $1$-Lipschitz with respect to $d$; namely,
for any $x$ and $y$ of ${\mathbb R}^m$, $| f(x) - f(y) | \leq d(x,y) $.
\end{Defi}
%
\begin{Defi} Let $r$ be a positive integer. For any points $x = (
x^{(1)},\ldots,\break x^{(2^r -1) } )'$ and $y = ( y^{(1)},\ldots, y^{(2^r -1) } )'$, we set
\[
d_{r} (x,y)= \sup_{j \in\{1,\ldots, 2^{r} -1 \}} \bigl| x^{(j)}
-y^{(j)} \bigr|.
\]
\end{Defi}
Let $L \in{\mathbb N}$ and $\ell\in\{1,\ldots, 2^{L-m(L)} \} $.
Let
\[
\Lambda_{Y,L} = \bigl(\Lambda_Y (s_j,
s_{j'} ) \bigr)_{ j,j'= 1,\ldots,
2^{r(L)} -1 },
\]
where the $ \Lambda_Y (s_j, s_{j'} )$ are defined in (\ref
{defgammaY}). Let
$G_{2^{m(L)} \Lambda_{Y,L}}$ denote the\break ${\mathcal N} ( 0, 2^{m(L)}
\Lambda_{Y,L})$-law and $P_{U_{L, \ell} | {\mathcal F}_{2^L + (\ell
-1)2^{m(L)}
}}$ be the conditional distribution of $U_{L,\ell}$ given $\mathcal
{F}_{2^L +
(\ell-1)2^{m(L)} }$.

According to \citet{Rus85} [see also Theorem 2 in
\citet{DedPriRay06}], there exists a random variable
$V_{L, \ell} = ( V^{(1)}_{L,\ell},\ldots,\break
V^{( 2^{r(L)} -1 )}_{L,\ell} )'$ with law $G_{2^{m(L)} \Lambda
_{Y,L}}$, measurable with respect to $\sigma(\delta_{2^L + \ell
2^{m(L)}}) \vee\sigma(U_{L,\ell}) \vee{\mathcal F}_{2^{L}+ (\ell
-1)2^{m(L)}} $, independent of ${\mathcal F}_{2^L + (\ell-1)2^{m(L)}}$
and such that
%
\begin{eqnarray}
\label{coupling1}
&&
{\mathbb E} \bigl( d_{r(L)} (U_{L,\ell},
V_{L,\ell}) \bigr) \nonumber\\
&&\qquad = {\mathbb E} \bigl( W_{d_{r(L)}}
(P_{U_{L, \ell} | {\mathcal F}_{2^L +
(\ell-1)2^{m(L)}}},G_{2^{m(L)} \Lambda_L}) \bigr)
\\
&&\qquad = \bkE\sup_{f \in\operatorname{Lip} (d_{r(L)} ) } \bigl( {\mathbb E} \bigl(
f(U_{L, \ell}) |
{{\mathcal F}_{2^L + (\ell-1)2^{m(L)}}} \bigr) - {\mathbb E} \bigl
(f(V_{L, \ell})
\bigr) \bigr).
\nonumber
\end{eqnarray}
By induction on $\ell$, the random variables $(V_{L,\ell})_{\ell= 1,\ldots, 2^{L-m(L)}}$
are mutually independent, independent of ${\mathcal F}_{2^L}$ and with
law ${\mathcal N} ( 0, 2^{m(L)}\Lambda_{Y,L})$. Hence, we have
constructed Gaussian random variables $(V_{L,\ell})_{L \in{\mathbb
N}, \ell=1,\ldots, 2^{L-m(L)}}$ that are mutually independent. In
addition, according to Lemma 2.11 of \citet{DudPhi83}, there
exists a Kiefer process $K_Y$ with covariance function $\Gamma_Y$ such
that for any $L \in{\mathbb N}$, any $\ell\in\{1,\ldots,
2^{L-m(L)} \} $ and any $j \in\{1,\ldots, 2^{r(L) - 1} \}$,
%
\begin{equation}
\label{relationVKiefer} V^{(j)}_{L, \ell} = K_Y
\bigl( s_j, 2^L + \ell2^{m(L)} \bigr)
-K_Y\bigl(s_j, 2^L + (\ell-1)
2^{m(L)}\bigr).
\end{equation}
Our construction is now complete.

In Proposition~\ref{propboundW1} proved in Section~\ref{sectionGApp},
we shall give some upper bounds for the quantities ${\mathbb E} (
d_{r(L)} (U_{L,\ell}, V_{L,\ell}) )$ for $L \in{\mathbb N}$
and $\ell\in\{1,\ldots,\break 2^{L-m(L)} \}$, showing that under our
condition on the dependence coefficients there exists a positive
constant $C$ such that
%
\begin{equation}
\label{evalerroraprox} {\mathbb E} \bigl( d_{r(L)} (U_{L,\ell},
V_{L,\ell}) \bigr) \leq C 2^{(m(L) +2r(L))/((2+\delta)\wedge3)} L^{2}.
\end{equation}
In Section~\ref{sectionerror} below, starting from (\ref{evalerroraprox}),
we bound up the error of approximation between the empirical process
and the Kiefer process.

\subsubsection{Upper bound for the approximation error} \label
{sectionerror} Let $\{ K_Y(s,t), s\in[0,1],\break t \geq0 \}$ be the
Gaussian process constructed as in step 1 with the following choice of
$r(L)$ and $m(L)$. For $\varepsilon<1/10$, let
%
\begin{equation}
\label{choicerm}\quad r(L) = \bigl([L/5]\wedge\bigl[ 2 \varepsilon L + 5
\log_2(L) \bigr]\bigr)\vee1 \quad\mbox{and}\quad m(L) = L - r(L),
\end{equation}
so that, for $L$ large enough,
%
\begin{eqnarray}
\label{encadrementrm}
2^{2 \varepsilon L - 1} L^{5} &\leq&2^{r(L)}
\leq2^{2 \varepsilon L} L^{5} \quad\mbox{and}\nonumber\\[-8pt]\\[-8pt]
2^{ L(1 - 2\varepsilon)} L^{-5}
&\leq&2^{m(L)} \leq2^{1+ L(1 - 2\varepsilon)} L^{- 5 }.\nonumber
\end{eqnarray}
Let $N\in{\mathbb N}^*$ and let $k \in{]1, 2^{N+1}]}$. To shorten the
notation, let $K_Y=K$ and \mbox{$R_Y=R$}. We first notice that
%
\begin{equation}
\label{1dec} \sup_{1 \leq k \leq2^{N+1} }\sup_{s \in
[0,1]} \bigl| R(s,k) - K(s,k) \bigr| \leq
\sup_{s \in[0,1]} \bigl| R(s,1) - K(s,1) \bigr| + \sum_{L =0}^N
D_L,\hspace*{-30pt}
\end{equation}
where
%
\begin{equation}
\label{defdl}\qquad\quad D_L:= \sup_{2^L < \ell\leq2^{L+1} }\sup_{s \in[0,1]} \bigl|
\bigl(R(s,\ell) -R\bigl(s,2^L\bigr) \bigr) - \bigl(K(s,\ell) -K
\bigl(s,2^L\bigr) \bigr) \bigr|.
\end{equation}
Notice first that $\sup_{s \in[0,1]} | R(s,1) - K(s,1) | \leq1 +\sup_{s
\in[0,1]} | K(s,1) | $. \citet{Ded10} (see the beginning of the
proof of his Theorem 2.1)
has proved that, for $u$ and $v$ in $ [0,1]$ and any positive integer $n$,
%
\begin{equation}
\label{rmkcov} \operatorname{Var} \bigl( K(u,n) -K (v, n ) \bigr) \leq C(\beta)
n |u-v|.
\end{equation}
Therefore, according to Theorem 11.17 in \citet{LedTal91},
$\bkE(\sup_{s \in[0,1]} | K(s,1) |) = O(1)$. It follows that for any
$\varepsilon\in{]0,1/2[}$,
%
\begin{equation}
\label{psd1} \sup_{s \in[0,1]} \bigl| R(s,1) - K(s,1) \bigr| = O\bigl( 2^{N (
1/2 -
\varepsilon)}
\bigr) \qquad\mbox{a.s.}
\end{equation}
To prove Theorem~\ref{mainemp}, it then suffices to prove that for any
$L \in\{0,\ldots, N \}$,
%
\begin{equation}
\label{butDl} D_L = O\bigl( 2^{L (1/2 - \varepsilon)}\bigr)
\qquad\mbox{a.s. for $\varepsilon= \delta^2/\bigl(22 (\delta+2)^2 \bigr)$.}
\end{equation}
With this aim, we decompose $D_L$ with the help of several quantities.
For any $K \in{\mathbb N}$ and any $s \in[0,1]$, let $\Pi_K (s) =
2^{-K} [2^K s ]$. Notice that the following decomposition is valid: for
any $L \in{\mathbb N}$,
%
\begin{equation}
\label{decsup} D_L \leq D_{L,1} + D_{L,2}+D_{L,3},
\end{equation}
where
\begin{eqnarray*}
D_{L,1} &=& \sup_{2^L < \ell\leq2^{L+1} }\sup_{s \in[0,1]} \bigl| \bigl
(R(s,\ell)
-R\bigl(\Pi_{r(L)}(s), \ell\bigr) \bigr)\\
&&\hspace*{73.4pt}{} - \bigl(R\bigl(s,2^L
\bigr) -R\bigl(\Pi_{r(L)}(s), 2^L \bigr) \bigr) \bigr|,
\\
D_{L,2} &=& \sup_{2^L < \ell\leq2^{L+1} }\sup_{s \in[0,1]} \bigl| \bigl
(K(s,\ell)
-K\bigl(\Pi_{r(L)}(s), \ell\bigr) \bigr) \\
&&\hspace*{73.4pt}{} - \bigl(K\bigl(s,2^L
\bigr) -K\bigl(\Pi_{r(L)}(s), 2^L \bigr) \bigr) \bigr|,
\\
D_{L,3} &=& \sup_{2^L < \ell\leq2^{L+1} }\sup_{s \in[0,1]} \bigl| \bigl
(R\bigl(
\Pi_{r(L)}(s),\ell\bigr)-R\bigl(\Pi_{r(L)}(s), 2^L
\bigr) \bigr)
\\
&&\hspace*{73.2pt}{} - \bigl(K\bigl(\Pi_{r(L)}(s),\ell\bigr) -K\bigl(\Pi_{r(L)}(s),
2^L\bigr) \bigr) \bigr|.
\end{eqnarray*}
In addition,
%
\begin{equation}
\label{decsup2} D_{L,3} \leq A_{L,3} + B_{L,3}
+C_{L,3},
\end{equation}
where
\begin{eqnarray*}
A_{L,3} &= &\sup_{j \in\{ 1,\ldots, 2^{r(L)-1} \}} \sup_{k \leq
2^{L-m(L)}} \Biggl| \sum
_{\ell=1}^k \bigl( U^{(j)}_{L,\ell} -
V^{(j)}_{L,\ell} \bigr) \Biggr|,
\\
B_{L,3} &=& \sup_{j \in\{ 1,\ldots, 2^{r(L)-1} \}} \sup_{k \leq
2^{L-m(L)}}
\sup_{ \ell\in I_{L,k}} \bigl| R(s_j, \ell) - R\bigl(s_j,
2^L + (k-1)2^{m(L)}\bigr) \bigr|,
\\
C_{L,3} &=&\sup_{j \in\{ 1,\ldots, 2^{r(L)-1} \}} \sup_{k \leq
2^{L-m(L)}}
\sup_{ \ell\in I_{L,k}} \bigl| K(s_j, \ell) - K\bigl(s_j,
2^L + (k-1)2^{m(L)}\bigr) \bigr|
\end{eqnarray*}
with $s_j = j2^{-r(L)}$.

Let us first deal with the terms $D_{L,2}$ and $C_{L,3}$ involving only
the approximating Kiefer process.
For any positive $\lambda$,
\begin{eqnarray*}
&&\p\bigl(| D_{L,2} | \geq\lambda\bigr)
\\
&&\qquad\leq\sum_{j =
1}^{2^{r(L)}} \p\Bigl(
\sup_{2^L < \ell\leq2^{L+1} } \sup_{s_{j-1}
\leq s \leq s_j} \bigl| \bigl(K(s,\ell) -K
\bigl(s,2^L\bigr) \bigr)\\
&&\hspace*{151pt}{} - \bigl(K(s_j, \ell) -K
\bigl(s_j, 2^L \bigr) \bigr) \bigr| \geq\lambda\Bigr).
\end{eqnarray*}
Setting
\begin{eqnarray*}
X(u,v)&=&
\bigl(K\bigl(s_j + u(s_{j+1} - s_j),2^L
+ v 2^L\bigr) -K\bigl(s_j + u(s_{j+1} -
s_j),2^L\bigr) \bigr) \\
&&{}- \bigl(K\bigl(s_j,
2^L + v 2^L \bigr) -K\bigl(s_j,
2^L \bigr) \bigr),
\end{eqnarray*}
we have
\[
\p( D_{L,2} \geq\lambda) \leq\sum_{j = 1}^{2^{r(L)}}
\p\Bigl( \sup_{(u,v)\in[0,1]^2 } \bigl| X(u,v) \bigr| \geq\lambda\Bigr).
\]
Using (\ref{rmkcov}), we infer that
\[
\bkE\bigl| X(u,v) -X\bigl(u',v'\bigr) \bigr|^2
\ll2^{L-r(L)} \bigl(\bigl| u-u'\bigr| + \bigl| v-v'\bigr| \bigr)
\]
and
\[
\sup_{(u,v)\in[0,1]^2}\bkE\bigl| X(u,v) \bigr|^2 \ll2^{L-r(L)}.
\]
Next, using Lemma 2 in \citet{Lai74}, as done in Lemma 6.2 in
\citet{BerPhi77}, and taking into account (\ref{encadrementrm}),
we infer that there exists a positive constant $c$ such that, for $L$
large enough,
\[
\p\bigl(|D_{L,2} | \geq c 2^{L(1/2 - \varepsilon)} \bigr) \ll
2^{r(L)} \exp\bigl( -L^{5}/2\bigr).
\]
Therefore,
%
\begin{equation}
\label{psd2} \sum_{L >0} \p\bigl( D_{L,2}
\geq c 2^{L(1/2 - \varepsilon)} \bigr) < \infty.
\end{equation}

Consider now the term $C_{L,3}$. For any positive $\lambda$,
\[
\p( C_{L,3} \geq\lambda) \leq\sum_{k =1}^{ 2^{L-m(L)}}
\p\Bigl( \sup_{s \in[0,1]} \sup_{
\ell\in I_{L,k}} \bigl| K(s, \ell) - K\bigl(s,
2^L + (k-1)2^{m(L)}\bigr) \bigr| \geq\lambda\Bigr).
\]
Setting $X(s,u)= K(s, 2^L + (k-1)2^{m(L)} +u 2^{m(L)}) - K(s, 2^L +
(k-1)2^{m(L)} +u 2^{m(L)})$ and using (\ref{rmkcov}), we have that
\[
\bkE\bigl| X(s,u) -X\bigl(s',u'\bigr) \bigr|^2
\ll2^{m(L)} \bigl( \bigl| s-s'\bigr| + \bigl| u-u'\bigr| \bigr)
\]
and
\[
\sup_{(s,u)\in[0,1]^2}\bkE\bigl| X(s,u) \bigr|^2 \ll2^{m(L)}.
\]
Therefore, by using once again Lemma 2 in \citet{Lai74}, as done in Lemma~6.3 in \citet{BerPhi77}, and taking into account (\ref
{encadrementrm}), we infer that there exists a positive constant $c$
such that, for $L$ large enough,
\[
\p\Bigl( \sup_{s \in[0,1]} \sup_{ \ell\in I_{L,k}} \bigl| K(s, \ell) -
K\bigl(s,
2^L + (k-1)2^{m(L)}\bigr) \bigr| \geq c 2^{L(1/2 - \varepsilon)} \Bigr)
\ll\exp\bigl(-L^{5}/2\bigr).
\]
Therefore,
%
\begin{equation}
\label{psd3} \sum_{L >0} \p\bigl( C_{L,3}
\geq c 2^{L(1/2 - \varepsilon)} \bigr) < \infty.
\end{equation}

We now prove that
%
\begin{equation}
\label{psd5} \sum_{L >0} \p\bigl( A_{L,3}
\geq2^{L(1/2 - \varepsilon)} \bigr) < \infty.
\end{equation}
From the stationarity of the sequence $((U_{L,\ell}, V_{L,\ell})
)_{\ell= 1,\ldots, 2^{L-m(L)}}$,
\[
\p\bigl( A_{L,3} \geq2^{L(1/2 - \varepsilon)} \bigr) \leq2^{L-m(L)}
2^{L( \varepsilon- 1/2)} {\mathbb E} \bigl( d_{r(L)} (U_{L,1},
V_{L,1}) \bigr).
\]
Therefore, by using (\ref{evalerroraprox}), we get that
\[
\p\bigl( A_{L,3} \geq2^{L(1/2 - \varepsilon)} \bigr) \ll2^{L(
\varepsilon- 1/2)}
2^{L-m(L)} 2^{{m(L) +2r(L)}/({(2+\delta)
\wedge3})} L^{2},
\]
which together with (\ref{encadrementrm}) proves (\ref{psd5}),
provided that
%
\begin{equation}
\label{firstconstraintepsi} \varepsilon<\frac{\delta\wedge1}{ 2 ( 8 + 3
( \delta\wedge1 ))}.
\end{equation}

We now show that
%
\begin{equation}
\label{psd4} \sum_{L >0} \p\bigl( B_{L,3}
\geq C 2^{L(1/2 - \varepsilon)} \bigr) < \infty.
\end{equation}
By stationarity, for any positive $\lambda$,
\[
\p( B_{L,3} \geq\lambda) \leq2^{L-m(L)} \sum
_{j =
1}^{2^{r(L)} } \p\Biggl(\sup_{ \ell\leq2^{m(L)} } \Biggl| \sum
_{i=1}^{\ell} \bigl( \I_{Y_i \leq j2^{-r(L)}} -
F_Y \bigl( j2^{-r(L)}\bigr) \bigr) \Biggr| \geq\lambda\Biggr).
\]
By Lemma~\ref{propY},
$|{\operatorname{Cov}} ( \I_{Y_0 \leq j2^{-r(L)}}, \I_{Y_i \leq
j2^{-r(L)}} )|
\leq{\mathbb E}(b(X_0, i)) = \beta(\sigma(X_0), X_i)$ and, consequently,
\[
\sum_{i \in{\mathbb Z}} \bigl| \operatorname{Cov} (
\I_{Y_0
\leq j2^{-r(L)}}, \I_{Y_i \leq j2^{-r(L)}} ) \bigr| \leq C(\beta).
\]
Applying Theorem 1 in \citet{DedMer10}, we get that for
any $v \geq1$,
\begin{eqnarray*}
&&\p\Biggl(\sup_{ \ell\leq2^{m(L)} } \Biggl| \sum_{i=1}^{\ell}
\biggl( \I_{Y_i \leq j2^{-r(L)}} - F_Y \biggl( \frac{j}{2^{r(L)}}\biggr)
\biggr) \Biggr| \geq4 \lambda\Biggr)
\\
&&\qquad\ll\biggl( 1 + \frac{\lambda^2}{2^{m(L)} v C(\beta)} \biggr)^{-v/4} +
\biggl( \frac{2^{m(L)}}{\lambda}
+ \frac{\lambda}{v } \biggr) \beta_{2,X} \biggl( \biggl[
\frac{\lambda}{v} \biggr] \biggr).
\end{eqnarray*}
Applying this inequality with $4\lambda= 2^{L(1/2 - \varepsilon)}$
and $v= L^{5}/C(\beta)$ and taking into account (\ref{encadrementrm})
together with our condition on the dependence coefficients, we derive
that for $L$ large enough,
\begin{eqnarray*}
&&
\p\Biggl(\sup_{ \ell\leq2^{m(L)} } \Biggl| \sum_{i=1}^{\ell}
\bigl( \I_{Y_i \leq j2^{-r(L)}} - F_Y \bigl( j2^{-r(L)}\bigr)
\bigr) \Biggr| \geq2^{L(
1/2 - \varepsilon)} \Biggr) \\
&&\qquad\ll\exp\bigl( -c_1
L^{5}\bigr) + L^{5\delta}2^{-L (1/2 - \varepsilon)\delta}.
\end{eqnarray*}
Therefore, (\ref{psd4}) holds provided that $
\varepsilon<\delta/(8 + 2\delta) $, which holds under~(\ref
{firstconstraintepsi}).

Taking into account (\ref{psd2}), (\ref{psd3}), (\ref{psd5}) and
(\ref{psd4}) together with the decompositions (\ref{decsup}) and
(\ref{decsup2}), the proof of (\ref{butDl}) will be complete if we
prove that,
for some positive constant $A$ to be chosen later,
%
\begin{equation}
\label{psd6} \sum_{L >0} \p\bigl( D_{L,1}
\geq\sqrt{AC(\beta)} 2^{L(1/2 - \varepsilon)} \bigr) < \infty.
\end{equation}
To shorten the notation, we set, for $\ell>m \geq0$,
%
\[
\mu_{\ell,m}(s) =R(s,\ell) - R(s,m)
\quad\mbox{and}\quad Z_{\ell, m}=d
\mu_{\ell,m}.
\]
We start from the elementary decomposition
\begin{eqnarray*}
&&
\mu_{\ell,2^L}(s)-\mu_{\ell,2^L}\bigl(\Pi_{r(L)}(s)\bigr)\\
&&\qquad=
\sum_{K=r(L)+1}^L \bigl(\mu_{\ell,2^L}\bigl(
\Pi_K(s)\bigr)-\mu_{\ell,2^L}\bigl(\Pi_{K-1}(s)\bigr)\bigr)
+ \mu_{\ell,2^L}(s)-\mu_{\ell,2^L}\bigl(\Pi_L(s)\bigr).
\end{eqnarray*}
Consequently,
%
\begin{equation}
\label{aj0}\quad \sup_{s \in[0, 1]} \bigl|\mu_{\ell,2^L}(s)-\mu_{\ell,2^L}
\bigl(\Pi_{r(L)}(s)\bigr)\bigr| \leq\sum_{K=r(L)+1}^L
\Delta_{K, \ell, 2^L} + \Delta_{L, \ell, 2^L}^*,
\end{equation}
where
\[
\Delta_{K, \ell, m}=\sup_{1 \leq i \leq2^K} \bigl|Z_{\ell,
m}
\bigl(\bigl](i-1)2^{-K}, i2^{-K}\bigr]\bigr)\bigr|
\]
and
\[
\Delta_{L, \ell,
m}^*=\sup_{s \in[0, 1]} \bigl|Z_{\ell, m}\bigl(\bigl]
\Pi_L(s), s\bigr]\bigr)\bigr|.\vadjust{\goodbreak}
\]
Note that
%
\begin{equation}
\label{aj1} - \bigl(\ell-2^{L}\bigr) {\mathbb P}\bigl(
\Pi_L(s)<Y_0\leq\Pi_L(s)+2^{-L}
\bigr) \leq Z_{\ell, 2^L}\bigl(\bigl]\Pi_L(s), s\bigr]\bigr)
\end{equation}
and
%
\begin{eqnarray}
\label{aj2} Z_{\ell, 2^L}\bigl(\bigl]\Pi_L(s), s\bigr]\bigr)&\leq&
Z_{\ell, 2^L}\bigl(\bigl]\Pi_L(s), \Pi_L(s)+2^{-L}\bigr]
\bigr)\nonumber\\[-8pt]\\[-8pt]
&&{} + \bigl(\ell-2^{L}\bigr) {\mathbb P}\bigl(\Pi_L(s)<Y_0
\leq\Pi_L(s)+2^{-L}\bigr).\nonumber
\end{eqnarray}
Applying Lemma~\ref{propY},
%
\begin{eqnarray}
\label{aj3}\qquad\quad {\mathbb P}\bigl(\Pi_L(s)<Y_0 \leq
\Pi_L(s)+2^{-L}\bigr)&\leq& C(\beta) {\mathbb
P}_0^*\bigl(\Pi_L(s)<Y_0 \leq
\Pi_L(s)+2^{-L}\bigr) \nonumber\\[-8pt]\\[-8pt]
&=& C(\beta) 2^{-L}.\nonumber
\end{eqnarray}
From (\ref{aj1}), (\ref{aj2}) and (\ref{aj3}), we infer that $\Delta
_{L, \ell, 2^L}^*\leq\Delta_{L, \ell, 2^L} +
C(\beta) $. Hence, it follows from (\ref{aj0}) that
\[
\sup_{s \in[0, 1]} \bigl|\mu_{\ell,2^L}(s)-\mu_{\ell,2^L}\bigl(
\Pi_{r(L)}(s)\bigr)\bigr| \leq C(\beta) + 2 \sum_{K=r(L)+1}^{L }
\Delta_{K, \ell, 2^L}.
\]
Therefore,
\begin{eqnarray*}
&&\sup_{2^L <\ell\leq2^{L+1}} \sup_{s \in[0, 1]} \bigl|\mu_{\ell,2^L}(s)-
\mu_{\ell,2^L}\bigl(\Pi_{r(L)}(s)\bigr)\bigr| \\
&&\qquad\leq C(\beta) + 2 \sum
_{K=r(L)+1}^{L} \sup_{2^L <\ell\leq2^{L+1}}
\Delta_{K, \ell, 2^L}.
\end{eqnarray*}
Hence, to prove (\ref{psd6}), it suffices to show that
%
\begin{equation}
\label{redux} \qquad\sum_{L>0} {\mathbb P} \Biggl( \sum
_{K=r(L)+1}^{L} \sup_{2^L <\ell
\leq2^{L+1}}
\Delta_{K, \ell, 2^L} > \sqrt{AC(\beta)}2^{L (
1/2 - \varepsilon)-2} \Biggr) < \infty.
\end{equation}
Let $c_K=(K(K+1))^{-1}$. Clearly, using the stationarity, (\ref
{redux}) is true provided that
%
\begin{equation}\qquad
\sum_{L >0} \sum_{K=r(L)+1}^{L}
{\mathbb P} \Bigl( \sup_{0 <\ell\leq
2^{L}}\Delta_{K, \ell, 0} >\sqrt{AC(\beta)}
c_K 2^{L ( 1/2 -
\varepsilon)-2} \Bigr) < \infty.
\end{equation}

We now give two upper bounds for the quantity
\[
{\mathbb P} \Bigl( \sup_{0 <\ell\leq2^{L}}\Delta_{K, \ell, 0} > \sqrt
{AC(\beta)}
c_K 2^{L ( 1/2 - \varepsilon)-2} \Bigr).
\]
Choose $p \in{]2,3]}$ such that $ p < 2(1+ \delta)$. Applying Markov's
inequality at order~$p$, we have
\[
{\mathbb P} \Bigl( \sup_{0 <\ell\leq2^{L}}\Delta_{K, \ell, 0} > \sqrt{AC(
\beta)}c_K 2^{L ( 1/2 - \varepsilon)-2} \Bigr) \ll c_K^{-p}
2^{L ( \varepsilon p -p/2 )} \Bigl\| \sup_{0 <\ell\leq
2^{L}} \Delta_{K, \ell, 0}
\Bigr\|_p^p.
\]
Applying inequality (7) of Proposition 1 in \citet{Wu07} to the
stationary sequence $(T_{K,i}^{ (j)})_{j \in{\mathbb Z}}$ defined by
$T_{K,i}^{ (j)} =\I_{(i-1)2^{-K}< Y_j \leq i2^{-K}}$, we have
\[
\Bigl\| \sup_{0 <\ell\leq2^{L}} \Delta_{K, \ell, 0} \Bigr\|_p
\leq2^{L/p} \sum_{j=0}^L
2^{-j/p} \| \Delta_{K, 2^j, 0} \|_p.
\]
Let $0<\eta< (p-2)/2$. \citet{Ded10} [see the displayed inequality
after (2.19) in his paper] proved that
\[
\label{boundDedec} \| \Delta_{K, 2^j, 0} \|_p^p
\ll2^{jp/2} \bigl( 2^{-K(p-2)/2} + 2^{-j\eta(2(1+\delta)-p)/2}+ 2^{j
\eta- j(p-2)/2
}
\bigr).
\]
Therefore,
%
\begin{eqnarray}
\label{firstbound}
&&\Bigl\| \sup_{0 <\ell\leq2^{L}}\Delta_{K, \ell, 0}
\Bigr\|_p^p \nonumber\\[-8pt]\\[-8pt]
&&\qquad\ll2^{L p/2 } \bigl( 2^{-K(p-2)/2} +
2^{-\eta L (2(1+\delta)-p)/2}+ 2^{\eta L -L(p-2)/2 } \bigr).\nonumber
\end{eqnarray}

On the other hand,
\begin{eqnarray*}
&&{\mathbb P} \Bigl( \sup_{0 <\ell\leq2^{L}}\Delta_{K, \ell, 0} > \sqrt
{AC(\beta)}
c_K 2^{L (1/2 - \varepsilon)-2} \Bigr)
\\
&&\qquad\leq\sum_{i=1}^{2^K} {\mathbb P} \Bigl(
\sup_{0 <\ell\leq2^{L}} \bigl|Z_{\ell,
0}\bigl(\bigl](i-1)2^{-K},
i2^{-K}\bigr]\bigr)\bigr|> \sqrt{AC(\beta)} c_K 2^{L (1/2 - \varepsilon)-2}
\Bigr).
\end{eqnarray*}
We now apply Theorem 1 in \citet{DedMer10}, taking into
account the stationarity: for any $x>0, v \geq1$, and $s_L^2 \geq2^L
\sum_{j=0}^{2^L} |{\operatorname{Cov}}(T_{K,i}^{ (0)}, T_{K,i}^{ (j)} )|$,
\begin{eqnarray*}
&&
{\mathbb P} \Bigl( \sup_{0 <\ell\leq2^{L}} \bigl|Z_{\ell,
0}\bigl(\bigl](i-1)2^{-K},
i2^{-K}\bigr]\bigr)\bigr|> 4x \Bigr) \\
&&\qquad\ll\biggl( \biggl( 1 + \frac
{x^2}{v s_L^2}
\biggr)^{-v/4} + 2^L \biggl( \frac1x + \frac{2x}{v
s_L^2}
\biggr) \beta_{2,X} \biggl( \biggl[\frac{x}{v} \biggr] \biggr)
\biggr).
\end{eqnarray*}
Applying Lemma~\ref{propY}, we have
$
|{\operatorname{Cov}}(T_{K,i}^{ (0)}, T_{K,i}^{ (j)} )|\leq2 \bkE(
T_{K,i}^{ (0)} b(X_0,j))
$. Hence,
%
\begin{equation}
\label{boundTKi} \sum_{j=0}^\infty\bigl|
\operatorname{Cov}\bigl(T_{K,i}^{ (0)}, T_{K,i}^{ (j)}
\bigr)\bigr| \leq C(\beta) {\mathbb P}_0^*\bigl((i-1)2^{-K}<
Y_0 \leq i2^{-K}\bigr) = C(\beta) 2^{-K}.\hspace*{-32pt}
\end{equation}
It follows that, for $K \geq r(L)$,
\[
\sum_{j=0}^\infty\bigl|\operatorname{Cov}
\bigl(T_{K,i}^{ (0)}, T_{K,i}^{ (j)} \bigr)\bigr|
\leq C(\beta) 2^{-r(L)}.
\]
For $L \geq2$, let $x=x_{K,L}= \sqrt{AC(\beta)} c_K 2^{L ( 1/2 -
\varepsilon)-4}$, $s_L^2 = C (\beta) 2^{L-r(L)} $ and $v=v_L=4L$.
Taking into account (\ref{encadrementrm}) and noting that $c_K \geq
(L(L+1))^{-1}$
for $K\leq L$, we obtain for $L$ large enough and $K\leq L$,
\[
\biggl(1 + \frac{x^2}{v s_L^2} \biggr)^{-v/4}\leq\biggl(1 +
\frac{A 2^{L(1-2\varepsilon)}}{2^{10} L^3(L+1)^2
2^{L-r(L)}} \biggr)^{-L} \leq3^{-L},
\]
the last bound being true provided $A$ is large enough. Hence,
for $L$ large enough and $ r(L) \leq K \leq L$,
%
\begin{eqnarray}
\label{secondbound}
&&{\mathbb P} \Bigl( \sup_{0 <\ell\leq2^{L}} \bigl|Z_{\ell,
0}
\bigl(\bigl](i-1)2^{-K}, i2^{-K}\bigr]\bigr)\bigr|> 4 x_{K,L}
\Bigr)\nonumber\\[-8pt]\\[-8pt]
&&\qquad\ll\biggl(\frac{1}{3^L} + \frac{L^{5+3\delta} 2^{L \varepsilon(2
+ \delta)}}{2^{L \delta/2}} \biggr).\nonumber
\end{eqnarray}

From (\ref{firstbound}) and (\ref{secondbound}), we then get that for
$L$ large enough
and any $\kappa\leq1$,
\begin{eqnarray*}
&&
\sum_{K=r(L)+1}^{L} {\mathbb P} \Bigl(
\sup_{0 <\ell\leq2^{L}}\Delta_{K, \ell, 0} > \sqrt{A C(\beta)} c_K
2^{L ( 1/2 - \varepsilon
)-2} \Bigr) \\
&&\qquad\ll\sum_{K=r(L)+1}^{[\kappa L]}
2^K \biggl( \frac{1}{3^{L}} + \frac{L^{5+3\delta} 2^{L \varepsilon(2 +
\delta)}}{2^{L \delta/2}} \biggr)
\\
&&\qquad\quad{}+ 2^{\varepsilon Lp} L^{2p} \sum_{K=[\kappa L]+1}^{L}
\bigl( 2^{-K(p-2)/2} + 2^{-\eta L (2(1+\delta)-p)/2}+ 2^{-L(p-2)/2+\eta L
} \bigr).
\end{eqnarray*}
Take $\kappa= \kappa(\varepsilon)= 1 \wedge2 \varepsilon(p+1) /(p-2)$.
It follows that (\ref{redux}) [and then (\ref{psd6})] holds provided
that the following constraints on $\varepsilon$ are satisfied:
\[
\varepsilon< \frac{p-2}{2(p+1)},\qquad \varepsilon\biggl(2 + \delta+
\frac{2(p+1)}{p-2} \biggr)< \delta/2,\qquad \varepsilon p < \frac
{p-2}{2} - \eta
\]
and
\[
\varepsilon p < \eta( 1+\delta-p/2).
\]
Let us take
\[
\eta= \frac{p-2}{4 + 2 \delta-p}
\quad\mbox{and}\quad p= 3 \wedge(2 + \delta/2).
\]
Both the above constraints on $\varepsilon$ and (\ref
{firstconstraintepsi}) are satisfied for $
\varepsilon= \delta^2/(22 (\delta+2)^2 ) $.
Therefore, (\ref{psd6}) holds, and Theorem~\ref{mainemp} follows.

\subsubsection{Gaussian approximation} \label{sectionGApp}
%
\begin{prop} \label{propboundW1}
For $L \in{\mathbb N}$, let $m(L) \in{\mathbb N}$ and $r(L) \in
{\mathbb N}^*$ be such that $ m (L) \leq L$ and $4r(L)\leq m(L)$.
Under the assumptions of Theorem~\ref{mainemp} and the notation of
Section~\ref{constructkiefer}, the following inequality holds: there
exists a positive constant $C$ not depending on $L$ such that, for any
$\ell\in\{ 1,\ldots, 2^{L-m(L)} \}$,
\[
{\mathbb E} \bigl( d_{r(L)} (U_{L,\ell}, V_{L,\ell})
\bigr) \leq C 2^{({m(L) +2r(L)})/({(2+\delta) \wedge3})} L^{2}.
\]
\end{prop}

\begin{pf}
From the stationarity of
the sequence $((U_{L,\ell}, V_{L,\ell}) )_{\ell= 1,\ldots,
2^{L-m(L)}}$, it suffices to prove the proposition for $\ell= 1$. Let
$L \in{\mathbb N}$ and $K \in\{0,\ldots,\break r(L) -1 \}$. To shorten the
notation, let us define the following set of integers:
\[
{\mathcal E}(L,K) = \bigl\{1,\ldots, 2^{r(L)-K}-1\bigr\} \cap(
2{\mathbb
N} +1 ),
\]
meaning that if $k \in{\mathcal E}(L,K) $, then $k$ is an odd integer
in $[1, 2^{r(L)-K}-1]$.

For $K \in\{0,\ldots, r(L) -1 \}$ and $k \in{\mathcal E}(L,K) $, define
\[
B_{K, k} = {\biggl]\frac{(k-1) 2^{K}}{2^{r(L)}}, \frac{k
2^{K}}{2^{r(L)}} \biggr]}
\quad\mbox{and}\quad Z^{(K,k)}_{L} = \sum_{i\in I_{L, 1}}
\bigl(\I_{Y_i \in B_{K,k}}-P_Y(B_{K,k}) \bigr).
\]
The associated column vector $Z_L$ in ${\mathbb R}^{2^{r(L)} -1}$ is
then defined by
\[
Z_{L } = \bigl( \bigl( Z^{(i,k_i)}_{L },
k_i\in{\mathcal E}(L,i) \bigr)_{i = 0,\ldots, r(L) -1 } \bigr)'.
\]
Notice that for any $j \in\{1,\ldots, 2^{r(L)} -1\}$,
%
\begin{equation}
\label{relationUZ} U^{(j)}_{L,1}= \sum
_{K=0}^{r(L)-1} \sum_{k_K \in{\mathcal E}(L,K)}
b_{K,k_K}(j)Z^{(K,k_K)}_{L }
\end{equation}
with $ b_{K,k_K}(j) = 0$ or $1$. This representation is unique in the
sense that, for $j$ fixed, there exists only one vector $ (
b_{(K,k_K)}(j), k_K\in{\mathcal E}(L,K) )_{K = 0,\ldots, r(L)
-1 } $ satisfying (\ref{relationUZ}). In addition, for any $K$ in $\{
0,\ldots, r(L) -1 \}$,
$\sum_{k\in{\mathcal E}(L,K)} b_{K,k}(j)\leq1$.
Let the column vector $b(j,L)$ and the matrix ${\mathbf P}_L $ be defined by
\[
b(j,L) = \bigl( \bigl( b_{K,k_K }(j), k_K \in{\mathcal
E}(L,K) \bigr)_{K = 0,\ldots, r(L) -1 } \bigr)'
\]
and
\[
{\mathbf P}_L = \bigl( b(1,L), b(2,L),\ldots, b
\bigl(2^{r(L)}-1, L\bigr) \bigr)'.
\]
${\mathbf P}_L$ has the following property: it is a square matrix of
${\mathbb R}^{2^{r(L)} -1}$ with determinant equal to $1$. Let us
denote by ${\mathbf P}_L^{-1}$ its inverse. With this notation, we then
notice that
%
\begin{equation}
\label{defZL} Z_L = {\mathbf P}_L^{-1}
U_{L,1}.
\end{equation}
Let now $a^2$ be a positive real and $V= ( V^{(1)}, \ldots,
V^{(2^{r(L)} - 1)} )'$ be a random variable with law ${\mathcal N}
( 0, a^2{\mathbf P}_L {\mathbf P}_L^T)$. According to the coupling relation
(\ref{coupling1}), we have that
%
\begin{eqnarray}
\label{maj1wasser}
{\mathbb E} \bigl( d_{r(L)} (U_{L,1},
V_{L,1}) \bigr) &=& {\mathbb E} \bigl( W_{d_{r(L)}}
(P_{U_{L, 1} | {\mathcal F}_{2^L }},G_{2^{m(L)}
\Lambda_L}) \bigr)
\nonumber
\\
&\leq& {\mathbb E} \bigl( W_{d_{r(L)}} (P_{U_{L, 1} |
{\mathcal F}_{2^L }} *
P_V,G_{2^{m(L)} \Lambda_L} * P_V) \bigr) \\
&&{}+ 2 {\mathbb E}
\bigl( d_{r(L)} (V,0) \bigr),\nonumber
\end{eqnarray}
where $*$ stands for the usual convolution product. Since $V^{(j)}$ is
a centered real Gaussian random variable with variance $v_j^2= a^2 \sum
_{K=0}^{r(L)-1} \sum_{k\in{\mathcal E}(L,K)} b_{K,k}(j)$, according
to inequality (3.6) in \citet{LedTal91}, we derive that
\begin{eqnarray*}
{\mathbb E} \bigl( d_{r(L)} (V,0) \bigr)&=& {\mathbb E} \Bigl(
\max_{j
\in\{ 1,\ldots,2^{r(L)} - 1 \} }\bigl|V^{(j)} \bigr| \Bigr)\\
&\leq&\bigl( 2 + 3
\bigl( \log
\bigl(2^{r(L)} - 1\bigr)\bigr)^{1/2} \bigr) \max_{j \in\{ 1,\ldots,2^{r(L)} -
1 \} }
v_j.
\end{eqnarray*}
Since $v_j^2 \leq a^2 r(L) \leq a^2 L $, we then get that
%
\begin{equation}
\label{boundVinfini} {\mathbb E} \bigl( d_{r(L)} (V,0) \bigr) \leq5 a L.
\end{equation}
Let us now give an upper bound for the quantity ${\mathbb E} (
W_{d_{r(L)}} (P_{U_{L, 1} | {\mathcal F}_{2^L }} * P_V,\break G_{2^{m(L)}
\Lambda_L} *P_V) )$ in (\ref{maj1wasser}). Let $(N_{i,L})_{i \in
{\mathbb Z}}$ be a sequence of independent random variables
with normal distribution ${\mathcal N}(0, \Lambda_L)$. Suppose,
furthermore, that the sequence $(N_{i,L})_{i \in{\mathbb Z}}$ is
independent of
${\mathcal F}_{\infty}\vee\sigma(\eta_i, i \in{\mathbb Z})$.
Denote by $\mathrm{I}_{2^{r(L)} -1}$ the identity matrix on ${\mathbb
R}^{2^{r(L)} -1}$ and let $N$ be a ${\mathcal N}(0,a^2
\mathrm{I}_{2^{r(L)} -1})$-distributed random variable, independent
of ${\mathcal F}_{\infty}\vee\sigma(N_{i,L}, i \in
{\mathbb Z}) \vee\sigma(\eta_i, i \in{\mathbb Z})$. Set $\widetilde
N_L = N_{1,L} + N_{2,L} + \cdots+ N_{2^{m(L)},L}$. We first notice that
%
\begin{eqnarray}
\label{boundWasser2first}\quad
& & {\mathbb E} \bigl( W_{d_{r(L)}}
(P_{U_{L, 1} | {\mathcal F}_{2^L }} * P_V,G_{2^{m(L)} \Lambda_L} *
P_V)
\bigr)
\nonumber\\[-8pt]\\[-8pt]
& &\qquad = \bkE\sup_{f \in\operatorname{Lip} (d_{r(L)} ) } \bigl( {\mathbb E} \bigl(
f(U_{L, 1} +{
\mathbf P}_L N ) | {{\mathcal F}_{2^L }} \bigr) - {\mathbb E}
\bigl(f(\widetilde N_{L} +{\mathbf P}_L N )\bigr) \bigr).\nonumber
\end{eqnarray}
Introduce now the following definition:
%
\begin{Defi} \label{defdstar} For two column vectors
\[
x = \bigl( \bigl(
x^{(i,k_i)}, k_i\in{\mathcal E}(L,i) \bigr)_{i =0,\ldots, r(L) -1 }
\bigr)'
\]
and $y = ( ( y^{(i,k_i)}, k_i\in{\mathcal E}(L,i) )_{i = 0,\ldots, r(L)
-1 } )'$ of ${\mathbb R}^{2^{r(L)} -1}$, let $d^*_{r(L)}$ be the
following distance:
\[
d^*_{r(L)} (x,y)= \sum_{K=0}^{r(L)-1}
\sup_{k \in{\mathcal E}(L,K)} \bigl| x^{(K,k)} -y^{(K,k)} \bigr|.
\]
Let also $\operatorname{Lip}(d^*_{r(L)}) $ be the set of functions from
${\mathbb R}^{2^{r(L)}-1}$ into ${\mathbb R}$ that are Lipschitz with
respect to $d^*_{r(L)} $, namely, $| f(x) - f(y) | \leq\sum
_{K=0}^{r(L)-1} \sup_{k \in{\mathcal E}(L,K)} | x^{(K,k)} -y^{(K,k)}
| $.
\end{Defi}

Let $x = ( x^{(1)}, \ldots, x^{(2^{r(L)} -1)} )'$ and $y =
( y^{(1)}, \ldots, y^{(2^{r(L)} -1)} )'$ be two column
vectors of ${\mathbb R}^{2^{r(L)} -1}$. Let now $u={\mathbf P}_L^{-1}x$ and
$v={\mathbf P}_L^{-1}y$. The vectors $u$ and $v$ of ${\mathbb R}^{2^{r(L)}
-1}$ can be rewritten $u = ( ( u^{(i,k_i)}, k_i\in{\mathcal
E}(L,i) )_{i =0,\ldots, r(L) -1 } )'$ and $v = (
( v^{(i,k_i)}, k_i\in{\mathcal E}(L,i) )_{i = 0,\ldots, r(L)
-1 } )' $. Notice now that if $f \in\operatorname{Lip}(d_{r(L)}) $, then
\begin{eqnarray*}
&&
\bigl|f(x) - f(y) \bigr| \\
&&\qquad\leq d_{r(L)} (x,y) = \sup_{j \in\{ 1,\ldots,2^{r(L)}
- 1 \}} \bigl|
b(j,L)' u - b(j,L)' v \bigr|
\\
&&\qquad\leq \sup_{j \in\{ 1,\ldots,2^{r(L)} - 1 \}} \sum_{K=0}^{r(L)-1}
\sum_{k_K \in{\mathcal E}(L,K) } b_{K,k_K}(j) \bigl| u^{(K,k_K)} -
v^{(K,k_K)}\bigr|
\\
&&\qquad\leq \sup_{j \in\{ 1,\ldots,2^{r(L)} - 1 \}} \sum_{K=0}^{r(L)-1}
\sum_{k_K \in{\mathcal E}(L,K) } b_{K,k_K}(j) \sup_{i
\in{\mathcal E}(L,K)}
\bigl| u^{(K,i)} - v^{(K,i)}\bigr|.
\end{eqnarray*}
Since for any $K \in\{0,\ldots, r(L) - 1 \}$ and any $j \in\{0,\ldots, 2^{r(L)} - 1 \}$,
\[
\sum_{k \in{\mathcal E}(L,K) } b_{K,k}(j)
\leq1,
\]
it follows that if $f \in\operatorname{Lip}(d_{r(L)}) $,
\begin{eqnarray*}
\bigl|f(x) - f(y) \bigr|&=&\bigl|f\circ{\mathbf P}_L (u)-f \circ{\mathbf P}_L
(v)\bigr| \leq\sum_{K=0}^{r(L)-1}
\sup_{k \in{\mathcal E}(L,K)} \bigl| u^{(K,k)} -v^{(K,k)} \bigr|\\
&=& d^*_{r(L)}
( u, v).
\end{eqnarray*}
Therefore, starting from (\ref{boundWasser2first}) and taking into
account (\ref{defZL}), we get
%
\begin{eqnarray}
\label{boundWasser2} & & {\mathbb E} \bigl( W_{d_{r(L)}} (P_{U_{L, 1} |
{\mathcal F}_{2^L }}
* P_V,G_{2^{m(L)} \Lambda_L} * P_V) \bigr)
\nonumber\\[-8pt]\\[-8pt]
& &\qquad \leq\bkE\sup_{f \in\operatorname{Lip}
(d^*_{r(L)} ) } \bigl( {\mathbb E} \bigl( f(Z_{L}
+N ) | {{\mathcal F}_{2^L }} \bigr) - {\mathbb E} \bigl(f\bigl({\mathbf
P}_L^{-1}\widetilde N_{L} + N\bigr)\bigr) \bigr).\nonumber
\end{eqnarray}
Let $\operatorname{Lip}(d^*_{r(L)}, {\mathcal F}_{2^L})$ be the set of
measurable functions $g\dvtx {\mathbb R}^{2^{r(L)} -1}
\times\Omega\rightarrow{\mathbb R}$ wrt the $\sigma$-fields
${\mathcal B} ( {\mathbb R}^{2^{r(L)} -1}) \otimes{\mathcal F}_{2^L} $
and ${\mathcal B} ({\mathbb R})$,
such that $g( \cdot, \omega) \in\operatorname{Lip}(d^*_{r(L)})$ and
$g(0,\omega)=0$ for any
$\omega\in\Omega$.
For the sake of brevity, we shall write $g(x)$ in place of $g(x, \omega
)$. From Point 2 of Theorem 1 in
\citet{DedPriRay06}, the following inequality holds:
%
\begin{eqnarray}
\label{equivalentdef}
&&\bkE\sup_{f \in\operatorname{Lip} (d^*_{r(L)} ) } \bigl(
{\mathbb E} \bigl(
f(Z_{L} +N ) | {{\mathcal F}_{2^L }} \bigr) - {\mathbb E}
\bigl(f\bigl({\mathbf P}_L^{-1}\widetilde N_{L} + N
\bigr)\bigr) \bigr)
\nonumber\\[-8pt]\\[-8pt]
&&\qquad= \sup_{g \in\operatorname{Lip}(d^*_{r(L)}, {\mathcal F}_{2^L})} {\mathbb E}
\bigl( g(Z_{L} +N) \bigr) - {
\mathbb E} \bigl(g\bigl({\mathbf P}_L^{-1}\widetilde
N_{L} + N\bigr)\bigr).
\nonumber
\end{eqnarray}
We shall prove that if $ a \in[L, L 2^{m(L)}]$, there exists a
positive constant $C$ not depending on $(L,a)$, such that
%
\begin{eqnarray}
\label{boundWasser3}
&&
\sup_{g \in\operatorname{Lip}(d^*_{r(L)}, {\mathcal
F}_{2^L})} {\mathbb E} \bigl( g(Z_{L}
+N) \bigr) - {\mathbb E} \bigl(g\bigl({\mathbf P}_L^{-1}
\widetilde N_{L} + N\bigr)\bigr) \nonumber\\
&&\qquad\leq C a^{-3}L^{5/2}2^{m(L)}
\nonumber\\[-8pt]\\[-8pt]
&&\qquad\quad{}+ C L^{-1} 2^{2r(L)} + C a^{-1-\delta} L^{\delta}
2^{2r(L) +m(L)} \nonumber\\
&&\qquad\quad{}+ C a^{-2 } L^2 2^{2r(L)+m(L)} +C
a^{-1} L^2 2^{r(L)}.
\nonumber
\end{eqnarray}
Gathering (\ref{boundWasser3}), (\ref{equivalentdef}), (\ref
{boundWasser2}), (\ref{maj1wasser}) and (\ref{boundVinfini}), and
taking
\[
a=L 2^{(m(L) +2r(L))/((2+\delta) \wedge3)},
\]
Proposition~\ref{propboundW1} will follow.

Let then $a\in[L, L 2^{m(L)}]$ and continue the proof by proving (\ref
{boundWasser3}). For any $i \geq1$, let $Y_{i,L}$ be the column vector
defined by $Y_{i,L}= ( Y^{(1)}_{i,L}, \ldots, Y^{(2^{r(L)}
-1)}_{i,L} )' $, where $Y^{(j)}_{i,L}= \I_{Y_{i +2^L} \leq
s_j}-F_Y(s_j)$. Notice then that
\[
Z_L = \sum_{i=1}^{2^{m(L)}}
Z_{i,L} \qquad\mbox{where }Z_{i,L} = {\mathbf P}_L^{-1}
Y_{i,L}.
\]
Therefore,
\[
Z_{i,L}= \bigl( \bigl( Z^{(K,k_K)}_{i,L}, k_K\in{\mathcal
E}(L,K) \bigr)_{K = 0,\ldots, r(L) -1 } \bigr)',
\]
where $Z^{(K,k)}_{i,L}= \I_{Y_{i +2^L} \in B_{K,k}}-P_Y(B_{K,k})$.
%
\begin{nota}\label{not21} Let $\varphi_a$ be the density of $N$ and
let for $x = ( ( x^{(i,k_i)}, k_i\in{\mathcal E}(L,K)
)_{i = 0,\ldots, r(L) -1 } )'$,
\[
g * \varphi_a (x,\omega) = \int g ( x+y, \omega)
\varphi_a(y) \,dy.
\]
For the sake of brevity, we shall write $g * \varphi_a (x)$
instead of $g * \varphi_a (x,\omega)$ (the partial derivatives will
be taken wrt $x$). Let also
\[
S_{0, L} = 0 \quad\mbox{and}\quad \mbox{for } j>0,\qquad S_{j, L} = \sum
_{i=1}^{j} Z_{i,L}.\vadjust{\goodbreak}
\]
\end{nota}
We now use the Lindeberg method to prove (\ref{boundWasser3}). We
first write that
%
\begin{eqnarray}
\label{lind1}
&&
\bkE\bigl( g (Z_L +N) - g \bigl({\mathbf
P}_L^{-1} \widetilde N_L +N \bigr) \bigr)
\nonumber\\
&&\qquad= \sum_{i=1}^{2^{m(L)}} \bkE\Biggl( g
\Biggl(S_{i-1, L} + Z_{i,L} + \sum_{j=i+1}^{2^{m(L)}}
{\mathbf P}_L^{-1} N_{j,L} + N \Biggr)\nonumber\\
&&\qquad\quad\hspace*{45pt}{} - g
\Biggl(S_{i-1, L} + {\mathbf P}_L^{-1} N_{i,L}
+ \sum_{j=i+1}^{2^{m(L)}} {\mathbf
P}_L^{-1} N_{j,L} + N \Biggr) \Biggr)
\\
&&\qquad\leq\sum_{i=1}^{2^{m(L)}} \sup_{g \in\operatorname{Lip}(d^*_{r(L)},
{\mathcal F}_{2^L})}
\bkE\bigl( g (S_{i-1, L} + Z_{i,L} + N ) \nonumber\\
&&\qquad\quad\hspace*{95pt}{}- g
\bigl(S_{i-1, L} + {\mathbf P}_L^{-1} N_{i,L} +
N \bigr) \bigr).
\nonumber
\end{eqnarray}
Let us introduce some notation and definitions.
%
\begin{Defi} \label{deftensorprod} For two positive integers $m$ and
$n$, let ${\mathcal M}_{m,n} ({\mathbb R})$ be the set of real matrices
with $m$ lines and $n$ columns. The Kronecker product
(or Tensor product) of $A = [ a_{i,j} ] \in{\mathcal M}_{m,n}
({\mathbb R})$ and $B = [ b_{i,j} ] \in{\mathcal M}_{p,q} ({\mathbb
R})$ is denoted by $A \otimes B$ and is defined to be the block matrix
\[
A \otimes B = \pmatrix{ a_{1,1}B & \cdots& a_{1,n}B
\cr
\vdots& & \vdots
\cr
a_{m,1}B & \cdots& a_{m,n}B} \in{\mathcal
M}_{mp,nq} ({\mathbb R}).
\]
For any positive integer $k$, the $k$th Kronecker power $A^{\otimes k}$
is defined inductively by $A^{\otimes1}=A$ and $A^{\otimes k} = A
\otimes A^{\otimes(k-1)}$.

If $\nabla$ denotes the differentiation operator given by $\nabla=
( \frac{\partial}{\partial x_1},\ldots, \frac{\partial
}{\partial x_m} )'$ acting on the differentiable functions $f\dvtx
{\mathbb R}^m \rightarrow{\mathbb R}$, we define
\[
\nabla\otimes\nabla= \biggl( \frac{\partial}{\partial x_1} \circ\nabla,\ldots,
\frac{\partial}{\partial x_m} \circ\nabla\biggr)'
\]
and $\nabla^{\otimes k}$ by $\nabla^{\otimes1}=\nabla$ and $\nabla
^{\otimes k} = \nabla\otimes\nabla^{\otimes(k-1)}$. If $f\dvtx
{\mathbb R}^m \rightarrow{\mathbb R}$ is $k$-times differentiable, for
any $x \in{\mathbb R}^m$, let $
D^k f(x) = \nabla^{\otimes k} f(x) $,
and for any vector $A$ of ${\mathbb R}^m$, we define $D^k f(x) \blddot
A^{\otimes k}$ as the usual scalar product in ${\mathbb R}^{m^k}$
between $D^k f(x)$ and $ A^{\otimes k}$.
\end{Defi}

For any $i\in\{1,\ldots, 2^{m(L)} \}$, let $G_{i,L} = {\mathbf P}_L^{-1}
N_{i,L}$,
\[
\Delta_{1,i,L} (g ) = g * \varphi_a (S_{i-1, L} +
Z_{i,L} ) - g * \varphi_a (S_{i-1, L} ) -
\tfrac{1}{2} D^2 g * \varphi_a (S_{i-1, L} )
\blddot G_{i,L}^{\otimes2}
\]
and
\[
\Delta_{2,i,L} (g) = g * \varphi_a (S_{i-1, L} +
G_{i,L} ) - g * \varphi_a (S_{i-1, L} ) -
\tfrac{1}{2} D^2 g * \varphi_a (S_{i-1, L} )
\blddot G_{i,L}^{\otimes2}.
\]
With this notation,
%
\begin{eqnarray}
\label{lind2}
&&\bkE\bigl( g (S_{i-1, L} + Z_{i,L} + N ) - g
\bigl(S_{i-1, L} + {\mathbf P}_L^{-1} N_{i,L} +
N \bigr) \bigr)\nonumber\\[-8pt]\\[-8pt]
&&\qquad= \bkE\bigl( \Delta_{1,i,L} (g) \bigr) - \bkE\bigl(
\Delta_{2,i,L} (g) \bigr).\nonumber
\end{eqnarray}
By the Taylor integral formula, noticing that $\bkE( G_{i,L}^{\otimes
3}) = 0$, we get
\[
\bigl| \bkE\bigl(\Delta_{2,i,L} (g) \bigr) \bigr| \leq\frac{1}{6} \biggl| \bkE
\int_0^1 D^4 g *
\varphi_a (S_{i-1, L} + t G_{i,L} ) \blddot
G_{i,L}^{\otimes4} \,dt \biggr|.
\]
Applying Lemma~\ref{lma4}, we then derive that
%
\begin{eqnarray}
\label{ajoutexpli1}
&&
\bigl| \bkE\bigl(\Delta_{2,i,L} (g) \bigr) \bigr| \nonumber\\
&&\qquad\ll
a^{-3} \bkE\Biggl( \Biggl( \sum_{K=0}^{r(L) -1}
\sup_{k\in{\mathcal E}(L,K)} \bigl|G_{1,L}^{(K,k)}\bigr| \Biggr) \Biggl( \sum
_{K=0}^{r(L) -1} \sum
_{k_K\in
{\mathcal E}(L,K)} \bigl(G_{1,L}^{(K,k_K)}
\bigr)^2 \Biggr)^{3/2} \Biggr)
\nonumber\\[-8pt]\\[-8pt]
&&\qquad\ll a^{-3} \Biggl( \bkE\Biggl( \sum_{K=0}^{r(L) -1}
\sup_{k\in
{\mathcal E}(L,K)} \bigl|G_{1,L}^{(K,k)}\bigr| \Biggr)^4
\Biggr)^{1/4} \nonumber\\
&&\qquad\quad{}\times\Biggl( \bkE\Biggl( \sum_{K=0}^{r(L) -1}
\sum_{k_K\in{\mathcal E}(L,K)} \bigl(G_{1,L}^{(K,k_K)}
\bigr)^2 \Biggr)^{2} \Biggr)^{3/4}.
\nonumber
\end{eqnarray}
Notice that
%
\begin{eqnarray}\quad
\label{ajoutexpli2} \sum_{K=0}^{r(L) -1}
\sup_{k\in{\mathcal E}(L,K)} \bigl|G_{1,L}^{(K,k)}\bigr| &\leq&\sum
_{K=0}^{r(L) -1} \biggl( \sum
_{k_K\in{\mathcal E}(L,K)} \bigl(G_{1,L}^{(K,k_K)}
\bigr)^2 \biggr)^{1/2}
\nonumber\\[-8pt]\\[-8pt]
&\leq&\sqrt{r(L)} \Biggl( \sum_{K=0}^{r(L) -1}
\sum_{k_K\in{\mathcal
E}(L,K)} \bigl(G_{1,L}^{(K,k_K)}
\bigr)^2 \Biggr)^{1/2}.
\nonumber
\end{eqnarray}
Moreover,
\begin{eqnarray*}
\bkE\Biggl( \sum_{K=0}^{r(L) -1} \sum
_{k_K\in{\mathcal E}(L,K)} \bigl(G_{1,L}^{(K,k_K)}
\bigr)^2 \Biggr)^{2} &\leq& \Biggl( \sum
_{K=0}^{r(L) -1} \sum_{k_K\in{\mathcal E}(L,K)}
\bigl( \bkE\bigl(G_{1,L}^{(K,k_K)}\bigr)^4
\bigr)^{1/2} \Biggr)^{2}
\\
&\leq& 3 \Biggl( \sum_{K=0}^{r(L) -1} \sum
_{k_K\in{\mathcal E}(L,K)} \bkE\bigl( \bigl(G_{1,L}^{(K,k_K)}
\bigr)^2\bigr) \Biggr)^{2}
\end{eqnarray*}
and
\[
\sum_{k\in{\mathcal E}(L,K)} \bkE\bigl( \bigl(G_{1,L}^{(K,k)}
\bigr)^2\bigr) = \sum_{k\in{\mathcal E}(L,K)} \biggl(\operatorname{Var}
\bigl( Z_{1,L}^{(K,k)}\bigr) +2 \sum _{i>0}\operatorname{Cov}
\bigl(Z_{1,L}^{(K,k)},Z_{i+1,L}^{(K,k)} \bigr) \biggr).
\]
Arguing as to get (\ref{boundTKi}), we then obtain that
\[
\sum_{k\in{\mathcal E}(L,K)} \bkE\bigl( \bigl(G_{1,L}^{(K,k)}
\bigr)^2\bigr) \leq C(\beta) \sum_{k\in{\mathcal E}(L,K)}
2^{K-r(L)} \leq C(\beta).
\]
From the above computations, it follows that
%
\begin{equation}
\label{ajoutexpli3} \bkE\Biggl( \sum_{K=0}^{r(L) -1}
\sum_{k_K\in{\mathcal E}(L,K)} \bigl(G_{1,L}^{(K,k_K)}
\bigr)^2 \Biggr)^{2} \leq3 \bigl( C(\beta)r(L)
\bigr)^2.
\end{equation}
Therefore, starting from (\ref{ajoutexpli1}), taking into account
(\ref{ajoutexpli2}), (\ref{ajoutexpli3}) and the fact that $r(L) \leq
L$, we then derive that
%
\begin{equation}
\label{lind3} \bigl| \bkE\bigl(\Delta_{2,i,L} (g) \bigr) \bigr| \ll
a^{-3} L^{5/2}.
\end{equation}
Let now
\begin{eqnarray*}
R_{1,i,L} (g ) &=& g * \varphi_a (S_{i-1, L} +
Z_{i,L} ) - g * \varphi_a (S_{i-1, L} ) - Dg *
\varphi_a (S_{i-1,
L} ) \blddot Z_{i,L} \\
&&{}-
\tfrac{1}{2} D^2 g * \varphi_a (S_{i-1, L} )
\blddot Z_{i,L}^{\otimes2}
\end{eqnarray*}
and
\begin{eqnarray*}
D_{1,i,L} (g) &=& Dg * \varphi_a (S_{i-1, L} ) \blddot
Z_{i,L} + \tfrac{1}{2} D^2 (g * \varphi_a)
(S_{i-1, L} ) \blddot Z_{i,L}^{\otimes2} \\
&&{}-
\tfrac{1}{2} D^2 g * \varphi_a (S_{i-1,
L} )
\blddot\bkE\bigl( G_{i,L}^{\otimes2} \bigr).
\end{eqnarray*}
With this notation,
%
\begin{equation}
\label{lind4} \bkE\bigl(\Delta_{1,i,L} (g) \bigr) = \bkE
\bigl(R_{1,i,L} (g) \bigr) + \bkE\bigl(D_{1,i,L} (g) \bigr).
\end{equation}
By the Taylor integral formula,
\[
\bigl|\bkE\bigl(R_{1,i,L} (g) \bigr) \bigr| \leq\biggl| \bkE\int_0^1
\frac
{(1-t)^2}{2} D^3 g * \varphi_a (S_{i-1, L}
+ t Z_{i,L} ) \blddot Z_{i,L}^{\otimes3} \biggr|.
\]
Applying Lemma~\ref{lma4} and using the fact that $\sup_{k\in
{\mathcal E}(L,K)} |Z_{i,L}^{(K,k)}| \leq2$ and $\sum_{k\in{\mathcal
E}(L,K)} (Z_{i,L}^{(K,k)})^2 \leq2$, we get that
%
\begin{equation}
\label{lind5first} \bigl| \bkE\bigl(R_{1,i,L} (g) \bigr) \bigr|\ll
a^{-2} \bigl(r(L)\bigr)^2 \ll a^{-2}
L^2.
\end{equation}
Let
%
\begin{equation}
\label{defgamma} \Delta(i,j) (g)= D^2 g * \varphi_a
(S_{i-j, L} ) - D^2 g * \varphi_a
(S_{i-j-1, L} )
\end{equation}
and
%
\begin{equation}
\label{defr} u_L = \bigl[a L^{-1}\bigr].
\end{equation}
Clearly, with the notation $X^{(0)} = X - \bkE(X)$,
%
\begin{eqnarray}
\label{lind5} D^2 g * \varphi_a (S_{i-1, L} )
\blddot\bigl( Z_{i,L}^{\otimes
2}\bigr)^{(0)}&=& \sum
_{j=1}^{(u_L\wedge i)-1} \Delta(i,j) (g) \blddot\bigl(
Z_{i,L}^{\otimes2}\bigr)^{(0)} \nonumber\\[-8pt]\\[-8pt]
&&{}+ D^2 g *
\varphi_a (S_{i-(u_L\wedge i), L} ) \blddot\bigl( Z_{i,L}^{\otimes
2}
\bigr)^{(0)}.\nonumber
\end{eqnarray}
For any $j\leq(u_L\wedge i)-1$, write
\[
\bkE\bigl( \Delta(i,j) (g)
\blddot\bigl( Z_{i,L}^{\otimes2}\bigr)^{(0)} \bigr) =
\bkE\bigl( \Delta(i,j) (g) \blddot\bkE_{i-j+2^L}\bigl(\bigl( Z_{i,L}^{\otimes
2}\bigr)^{(0)}\bigr) \bigr)
\]
and notice that, by Lemma~\ref{lma5},
\begin{eqnarray*}
&&
\bkE\bigl( \Delta(i,j) (g) \blddot\bkE_{i-j+2^L}\bigl(
Z_{i,L}^{\otimes
2}\bigr)^{(0)} \bigr)
\\
&&\qquad\leq\sup_{t \in[0,1]} \bigl| \bkE\bigl( D^3 g *
\varphi_a (S_{i-j-1, L} +tZ_{i-j,L}) \blddot\bigl(
Z_{i-j,L}\otimes\bkE_{i-j+2^L}\bigl( Z_{i,L}^{\otimes2}
\bigr)^{(0)} \bigr) \bigr) \bigr|
\\
&&\qquad\ll a^{-2} \sum_{K_1,k_{K_1}} \sum
_{K_2,k_{K_2}}\sum_{K_3,k_{K_3}} \bkE
\bigl(\bigl|Z_{i-j,L}^{K_1,k_{K_1}}\bigr| \bigl| \bkE_{i-j+2^L}\bigl(
Z_{i,L}^{K_2,k_{K_2}} Z_{i,L}^{K_3,k_{K_3}} \\
&&\hspace*{235.6pt}{}- \bkE
\bigl(Z_{i,L}^{K_2,k_{K_2}} Z_{i,L}^{K_3,k_{K_3}} \bigr)
\bigr) \bigr| \bigr),
\end{eqnarray*}
where $K_i \in\{ 0,\ldots, r(L) -1 \}$ and $k_{K_i} \in{\mathcal
E}(L,K_i)$, for any $i\in\{1,2,3\}$. Applying Lemma~\ref{propY}, we
infer that
\[
\bigl| \bkE_{i-j+2^L}\bigl( Z_{i,L}^{K_2,k_{K_2}}
Z_{i,L}^{K_3,k_{K_3}} - \bkE\bigl(Z_{i,L}^{K_2,k_{K_2}}
Z_{i,L}^{K_3,k_{K_3}} \bigr)\bigr) \bigr| \leq4 b_1 \bigl( {
\mathcal F}_{i-j+ 2^L}, i + 2^L \bigr).
\]
Since $\sum_{K_1=0}^{r(L)-1}\sum_{k_{K_1} \in{\mathcal E}(L,K_1)}
|Z_{i-j,L}^{K_1,k_{K_1}}| \leq2 r(L)$ and $\bkE( b_1( {\mathcal
F}_{i-j+ 2^L}, i + 2^L ) ) \leq\beta_{1,X} (j) $, we then derive that
%
\begin{equation}
\label{lind6} \bkE\bigl( \Delta(i,j) (g) \blddot\bigl( Z_{i,L}^{\otimes2}
\bigr)^{(0)} \bigr) \ll a^{-2} r(L) 2^{2r(L)}
\beta_{1,X} (j).
\end{equation}
On the other hand, by using Lemma~\ref{lma5}, we infer that
\begin{eqnarray*}
&&
\bkE\bigl( D^2 g * \varphi_a (S_{i-(u_L\wedge i), L} )
\blddot\bigl( Z_{i,L}^{\otimes2}\bigr)^{(0)} \bigr) \\
&&\qquad= \bkE
\bigl( D^2 g * \varphi_a (S_{i-(u_L\wedge i), L} ) \blddot
\bkE_{i-(u_L
\wedge i) + 2^L}\bigl( Z_{i,L}^{\otimes2}\bigr)^{(0)}
\bigr)
\\
&&\qquad\ll a^{-1} \sum_{K_1,k_{K_1}} \sum
_{K_2,k_{K_2}} \bkE\bigl( \bigl| \bkE_{i-(u_L \wedge i) + 2^L}\bigl(
Z_{i,L}^{K_1,k_{K_1}} Z_{i,L}^{K_1,k_{K_1}} \\
&&\hspace*{185pt}{} - \bkE
\bigl(Z_{i,L}^{K_1,k_{K_1}} Z_{i,L}^{K_1,k_{K_1}} \bigr)
\bigr) \bigr| \bigr).
\end{eqnarray*}
Using the same arguments as to get (\ref{lind6}), we obtain that
%
\begin{equation}
\label{lind7} \bkE\bigl( D^2 g * \varphi_a
(S_{i-(u_L\wedge i), L} ) \blddot\bigl( Z_{i,L}^{\otimes2}
\bigr)^{(0)} \bigr)\ll a^{-1} 2^{2r(L)}
\beta_{1,X} (u_L \wedge i).
\end{equation}
Starting from (\ref{lind5}) and taking into account (\ref{lind6}),
(\ref{lind7}), the choice of $u_L$ and the condition on the $\beta
$-dependence coefficients, we then derive that
%
\begin{eqnarray}
\label{lind8}
&&\sum_{i=1}^{2^{m(L)}} \bkE\bigl(
D^2 g * \varphi_a (S_{i-1,
L} ) \blddot\bigl(
Z_{i,L}^{\otimes2}\bigr)^{(0)} \bigr)\nonumber\\[-8pt]\\[-8pt]
&&\qquad\ll2^{2r(L)} a^{-1} \biggl( \frac{2^{m(L)} L^{1 + \delta}}{a^{1 + \delta
}} +
2^{m(L)} \frac{L}{a} \biggr).\nonumber
\end{eqnarray}
To give now an estimate of the expectation of $Dg * \varphi_a
(S_{i-1, L} ) \blddot Z_{i,L} $, we write
\[
Dg * \varphi_a (S_{i-1, L} )= Dg * \varphi_a (0) + \sum_{j=1}^{i-1}
\bigl( Dg *
\varphi_a (S_{i-j, L} )- Dg * \varphi_a
(S_{i-j-1, L}) \bigr).
\]
Hence,
%
\begin{eqnarray}
\label{lind9}
&&\bkE\bigl(Dg * \varphi_a (S_{i-1, L} ) \blddot
Z_{i,L} \bigr) \nonumber\\
&&\qquad= \bkE\bigl( Dg * \varphi_a (0 ) \blddot
Z_{i,L} \bigr) \\
&&\qquad\quad{}+ \sum_{j=1}^{i-1}
\bkE\bigl( \bigl( Dg * \varphi_a (S_{i-j, L})- Dg *
\varphi_a (S_{i-j-1, L} ) \bigr)\blddot Z_{i,L} \bigr).\nonumber
\end{eqnarray}
Applying Lemma~\ref{propY},
\begin{eqnarray*}
\bigl|\bkE\bigl( Dg * \varphi_a (0 ) \blddot Z_{i,L} \bigr)\bigr| &=&
\bigl|\bkE\bigl( Dg * \varphi_a (0 ) \blddot\bkE_{ 2^L} (
Z_{i,L} )\bigr)\bigr|
\\
&\leq&\bkE\Biggl( \sum_{K =0}^{r(L) -1} \sum
_{k_K\in{\mathcal
E}(L,K)} \biggl| \frac{\partial g * \varphi_a }{\partial x^{(K,k_K)}} (0) \biggr| b_1
\bigl( {\mathcal F}_{ 2^L}, i + 2^L \bigr) \Biggr).
\end{eqnarray*}
Notice now that by inequality (\ref{inelma1}), for any $K $ in $\{
0,\ldots, r(L)-1 \}$,
the random variable
\[
\sum_{k\in{\mathcal E}(L,K)} \biggl| \frac{\partial g * \varphi_a
}{\partial x^{(K,k)}} (0) \biggr|
\]
is a ${\mathcal F}_{2^L}$-measurable random variable with infinite norm
less than one. Therefore,
%
\begin{equation}
\label{lind9bis} \bigl|\bkE\bigl( Dg * \varphi_a (0 ) \blddot
Z_{i,L} \bigr)\bigr| \ll r(L) \beta_{1,X} (i).
\end{equation}
We give now an estimate of $\sum_{j=1}^{i-1} \bkE( ( Dg *
\varphi_a (S_{i-j, L})- Dg * \varphi_a (S_{i-j-1, L} ) ) \blddot
Z_{i,L} ) $. By Lemmas~\ref{lma5} and~\ref{propY}, for any
$i \geq j+1$,
\begin{eqnarray*}
\hspace*{-5.5pt}&&
\bigl|\bkE\bigl( \bigl( Dg * \varphi_a (S_{i-j, L})- Dg *
\varphi_a (S_{i-j-1, L} ) \bigr) \blddot Z_{i,L} \bigr)
\bigr|
\\
\hspace*{-5.5pt}&&\quad= \bigl|\bkE\bigl( \bigl( Dg * \varphi_a (S_{i-j, L})- Dg *
\varphi_a (S_{i-j-1, L} ) \bigr) \blddot\bkE_{i-j + 2^L}
(Z_{i,L} ) \bigr)\bigr|
\\
\hspace*{-5.5pt}&&\quad\leq\sup_{t \in[0,1]} \bigl| \bkE\bigl( D^2 g *
\varphi_a (S_{i-j-1, L} +tZ_{i,L}) \blddot\bigl(
Z_{i-j,L}\otimes\bkE_{i-j+2^L}( Z_{i,L})\bigr) \bigr) \bigr|
\\
\hspace*{-5.5pt}&&\quad\ll a^{-1} \sum_{K_1=0}^{r(L)-1}\sum
_{k_{K_1} \in{\mathcal
E}(L,K_1)} \sum_{K_2=0}^{r(L)-1}
\sum_{k_{K_2} \in{\mathcal E}(L,K_2)} \bkE\bigl(\bigl|Z_{i-j,L}^{K_1,k_{K_1}}\bigr|
b_1\bigl( {\mathcal F}_{i-j+ 2^L}, i + 2^L \bigr)
\bigr).
\end{eqnarray*}
We then infer that for any $i \geq j+1$,
%
\begin{eqnarray}
\label{lind9ter}
&&\bigl| \bkE\bigl( \bigl( Dg * \varphi_a
(S_{i-j, L})- Dg * \varphi_a (S_{i-j-1, L} ) \bigr)
\blddot Z_{i,L} \bigr) \bigr| \nonumber\\[-8pt]\\[-8pt]
&&\qquad\ll a^{-1} r(L) 2^{r(L)}
\beta_{1,X} (j).\nonumber
\end{eqnarray}
From now on, we assume that $j<i \wedge u_L $. Notice that
\begin{eqnarray*}
&&
\bigl( Dg * \varphi_a (S_{i-j, L})- Dg *
\varphi_a (S_{i-j-1, L} ) \bigr) \blddot Z_{i,L}\\
&&\quad{} -
D^2g * \varphi_a (S_{i-j-1, L} )\blddot
(Z_{i-j,L} \otimes Z_{i,L} )
\\
&&\qquad= \int_0^1 (1- t) D^3 g *
\varphi_a (S_{i-j-1, L} +t Z_{i-j,L}) \blddot
\bigl(Z_{i-j,L}^{\otimes2} \otimes Z_{i,L} \bigr) \,dt.
\end{eqnarray*}
By using Lemmas~\ref{lma5} and~\ref{propY}, we infer that
\begin{eqnarray*}
&&
\biggl| \bkE\biggl(\int_0^1 (1- t)
D^3 g * \varphi_a (S_{i-j-1, L} +t
Z_{i-j,L}) \blddot\bigl(Z_{i-j,L}^{\otimes2} \otimes
Z_{i,L} \bigr) \,dt \biggr) \biggr|
\\
&&\qquad\ll a^{-2} \sum_{K_1=0}^{r(L)-1}\sum
_{k_{K_1} \in{\mathcal E}(L,K_1)} \sum_{K_2=0}^{r(L)-1}
\sum_{k_{K_2} \in{\mathcal E}(L,K_2)} \sum_{K_3=0}^{r(L)-1}
\sum_{k_{K_3} \in{\mathcal E}(L,K_3)}
\\
&&\qquad\quad\hspace*{130pt}\bkE\bigl(\bigl|Z_{i-j,L}^{K_1,k_{K_1}}\bigr| \bigl|Z_{i-j,L}^{K_2,k_{K_2}}\bigr|b_1
\bigl( {\mathcal F}_{i-j+ 2^L}, i + 2^L \bigr) \bigr).
\end{eqnarray*}
Therefore,
%
\begin{eqnarray}
\label{lind10}
&&\biggl| \bkE\biggl(\int_0^1 (1- t)
D^3 g * \varphi_a (S_{i-j-1, L} +t
Z_{i-j,L}) \blddot\bigl(Z_{i-j,L}^{\otimes2} \otimes
Z_{i,L} \bigr) \,dt \biggr) \biggr| \nonumber\\[-8pt]\\[-8pt]
&&\qquad\ll a^{-2} \bigl(r(L)
\bigr)^2 2^{r(L)} \beta_{1,X} (j).\nonumber
\end{eqnarray}
In order to estimate the term $ \bkE( D^2 g * \varphi_a
(S_{i-j-1, L}) \blddot( Z_{i-j,L} \otimes Z_{i,L}) )$,
we use the following decomposition:
\begin{eqnarray*}
&&
D^2 g * \varphi_a (S_{i-j-1, L}) \\
&&\qquad= \sum
_{l=1}^{(j-1) \wedge(i-j-1)} \bigl( D^2 g *
\varphi_a (S_{i-j-l,L}) - D^2 g *
\varphi_a (S_{i-j-l-1,L}) \bigr) \\
&&\qquad\quad{}+ D^2 g *
\varphi_a (S_{(i-2j)\vee0, L}).
\end{eqnarray*}
For any $l\in\{ 1, \ldots, (j-1) \wedge(i-j-1) \}$, using the same
arguments as to get (\ref{lind10}), we obtain that
%
\begin{eqnarray}
\label{lind11}\qquad
&&\bigl|\bkE\bigl( \bigl( D^2 g * \varphi_a
(S_{i-j-l,L}) - D^2 g * \varphi_a
(S_{i-j-l-1,L}) \bigr) \blddot(Z_{i-j,L} \otimes Z_{i,L} )
\bigr)\bigr| \nonumber\\[-8pt]\\[-8pt]
&&\qquad\ll a^{-2} \bigl(r(L)\bigr)^2 2^{r(L)}
\beta_{1,X} (j).\nonumber
\end{eqnarray}

As a second step, we bound up
$ | \bkE( D^2 g * \varphi_a (S_{(i-2j)\vee0,L}) \blddot
(Z_{i-j,L} \otimes Z_{i,L} )^{(0)} ) |$. Assume first that $j
\leq[i/2]$. Clearly, using the notation (\ref{defgamma}),
\[
D^2 g * \varphi_a (S_{i-2j,L})= \sum
_{l=j}^{(u_L-1) \wedge(i-j -1)} \Delta(i, l+j) (g)+ D^2 g *
\varphi_a (S_{(i-j-u_L) \vee0,L}).
\]
Now for any $l\in\{ j,\ldots, (u_L-1) \wedge(i-j-1) \}$, by using
Lemma~\ref{lma5}, we get that
\begin{eqnarray*}
&&
\bigl| \bkE\bigl( \Delta(i, l+j) \blddot(Z_{i-j,L} \otimes
Z_{i,L} )^{(0)} \bigr)\bigr|
\\
&&\qquad
\ll a^{-2} \sum_{K_1,k_{K_1} } \sum
_{K_2,k_{K_2} } \sum_{K_3,k_{K_3}} \bkE\bigl|
Z_{i-j-l,L}^{K_1,k_{K_1}} \bkE_{i-j - l + 2^L}\bigl( Z_{i-j,L}^{K_2,k_{K_2}}
Z_{i,L}^{K_3,k_{K_3}} \\
&&\qquad\quad\hspace*{205pt}{} - \bkE\bigl(Z_{i-j,L}^{K_2,k_{K_2}}
Z_{i,L}^{K_3,k_{K_3}} \bigr)\bigr) \bigr|.
\end{eqnarray*}
Applying Lemma~\ref{propY}, we infer that
\begin{eqnarray*}
&&\bigl| \bkE_{i-j - l + 2^L}\bigl( Z_{i-j,L}^{K_2,k_{K_2}}
Z_{i,L}^{K_3,k_{K_3}} - \bkE\bigl(Z_{i-j,L}^{K_2,k_{K_2}}
Z_{i,L}^{K_3,k_{K_3}} \bigr)\bigr) \bigr| \\
&&\qquad\leq4 b_2 \bigl( {
\mathcal F}_{i-j - l +
2^L}, i-j+2^L, i+2^L \bigr).
\end{eqnarray*}
Therefore,
%
\begin{equation}
\label{lind12} \bigl| \bkE\bigl( \Delta(i, l+j) \blddot(Z_{i-j,L} \otimes
Z_{i,L} )^{(0)} \bigr)\bigr| \ll a^{-2} r(L)
2^{2r(L)} \beta_{2,X} (l).
\end{equation}
If $j \leq i-u_L$, with similar arguments,
%
\begin{equation}
\label{lind13} \qquad\bigl| \bkE\bigl( D^2 g * \varphi_a
(S_{i-j-u_L,L } ) \blddot(Z_{i-j,L} \otimes Z_{i,L}
)^{(0)} \bigr)\bigr| \ll a^{-1}2^{2r(L)}
\beta_{2,X} (u_L).
\end{equation}
Now if $j > i-u_L$, we infer that
%
\begin{equation}
\label{lind14} \bigl| \bkE\bigl(\bigl( D^2 g * \varphi_a (0 )
\bigr) \blddot(Z_{i-j,L} \otimes Z_{i,L} )^{(0)} \bigr)\bigr|
\ll a^{-1}2^{2r(L)} \beta_{2,X} \bigl([i/2]\bigr)
\end{equation}
by using also the fact that, since $j\leq[i/2]$, $\beta_{2,X} (i-j)
\leq\beta_{2,X} ([i/2]) $.
Assume now that $j \geq[i/2] + 1$. For any $j \leq i$, we get
%
\begin{equation}
\label{lind15} \bigl| \bkE\bigl(\bigl( D^2 g * \varphi_a (0 )
\bigr) \blddot Z_{i-j,L} \otimes Z_{i,L} \bigr)\bigr| \ll a^{-1}
r(L) 2^{r(L)} \beta_{1,X} \bigl([i/2]\bigr).
\end{equation}
Starting from (\ref{lind9}), adding inequalities (\ref
{lind9bis})--(\ref{lind15}) and summing on $j$ and~$l$,
we then obtain
\begin{eqnarray*}
&&
\Biggl|\bkE\bigl(Dg * \varphi_a (S_{i-1, L} ) \blddot
Z_{i,L} \bigr)\\
&&\quad{} - \sum_{j=1}^{u_L-1}
\bkE\bigl(D^2 g * \varphi_a (S_{i-2j, L})\bigr)
\blddot\bkE(Z_{i-j,L} \otimes Z_{i,L} ) \I_{j \leq[i/2]} \Biggr|
\\
&&\qquad
\ll r(L) \beta_{1,X}(i) + a^{-1}L 2^{r(L)} \sum
_{j=u_L}^{i} \beta_{1,X} (j) +
a^{-1} 2^{2 r(L)} u_L \beta_{2,X}(u_L)
\\
&&\qquad\quad{}+ a^{-1}2^{2r(L)} u_L \beta_{2,X}
\bigl([i/2]\bigr) + a^{-2} L 2^{2r(L)} \sum
_{j=1}^{ u_L} j \beta_{2,X} (j).
\end{eqnarray*}
Next, summing on $i$ and taking into account the condition on the
$\beta$-dependence coefficients and the choice of $u_L$, we get that
%
\begin{eqnarray}
\label{lind16}
&&\sum_{i =1}^{2^{m(L)}} \Biggl|\bkE
\bigl(Dg * \varphi_a (S_{i-1, L} ) \blddot Z_{i,L}
\bigr) \nonumber\\
&&\hspace*{10pt}\quad{}- \sum_{j=1}^{u_L-1} \bkE
\bigl(D^2 g * \varphi_a (S_{i-2j, L})\bigr) \blddot
\bkE(Z_{i-j,L} \otimes Z_{i,L} ) \I_{j \leq[i/2]} \Biggr|
\\
&&\hspace*{10pt}\qquad\ll L^{-1} 2^{2r(L)} + a^{-1-\delta} L^{\delta}
2^{2r(L) +m(L)} + a^{-2 } L^2 2^{2r(L)+m(L)}.
\nonumber
\end{eqnarray}
It remains to bound up
\begin{eqnarray*}
A_i&:=& \Biggl|\sum_{j=1}^{u_L-1} \bkE
\bigl(D^2 g * \varphi_a (S_{i-2j}) \bigr) \blddot
\bkE(Z_{i-j,L} \otimes Z_{i,L} ) \I_{j
\leq[i/2]} \\
&&\hspace*{34pt}{}- \sum
_{j=1}^{\infty} \bkE\bigl(D^2 g *
\varphi_a (S_{i-1}) \bigr)\blddot\bkE(Z_{i-j,L}
\otimes Z_{i,L} ) \Biggr|.
\end{eqnarray*}
We first notice that by Lemma~\ref{lma5}, for any positive integer $j$,
\begin{eqnarray*}
&&\bigl|\bkE\bigl(D^2 g * \varphi_a (S_{i-1})\bigr)
\blddot\bkE(Z_{i-j,L} \otimes Z_{i,L} )\bigr|
\\
&&\qquad\ll a^{-1} \sum_{K_1=0}^{r(L)-1}\sum
_{k_{K_1} \in{\mathcal E}(L,K_1)} \sum_{K_2=0}^{r(L)-1}
\sum_{k_{K_2} \in{\mathcal E}(L,K_2)} \bigl| \bkE\bigl(Z_{i-j,L}^{K_1,k_{K_1}}
\bkE_{i-j+2^L} \bigl(Z_{i,L}^{K_2,k_{K_2}} \bigr) \bigr) \bigr|.
\end{eqnarray*}
Therefore,
%
\begin{equation}
\label{lind17}\quad \bigl|\bkE\bigl(D^2 g * \varphi_a
(S_{i-1})\bigr)\blddot\bkE(Z_{i-j,L} \otimes
Z_{i,L})\bigr| \ll a^{-1} r(L) 2^{r(L)}
\beta_{1,X}(j).
\end{equation}
On an other hand, applying Lemma~\ref{lma5}, we obtain for any $i \geq
2$ and any $j \in\{1,\ldots, [i/2] \}$,
\begin{eqnarray*}
\hspace*{-3pt}&&
\bigl|\bkE\bigl( D^2 g * \varphi_a (S_{i-1}) -
D^2 g * \varphi_a (S_{i-2j})\bigr) \blddot\bkE
(Z_{i-j,L} \otimes Z_{i,L} )\bigr|
\\
\hspace*{-3pt}&&\qquad\ll a^{-2} \sum_{K_1=0}^{r(L)-1}\sum
_{k_{K_1} \in{\mathcal E}(L,K_1)} \sum_{K_2=0}^{r(L)-1}
\sum_{k_{K_2} \in{\mathcal E}(L,K_2)} \sum_{K_3=0}^{r(L)-1}
\sum_{k_{K_3} \in{\mathcal E}(L,K_3)} \sum_{\ell
=1}^{2j-1}
\\
\hspace*{-3pt}&&\qquad\quad\hspace*{130pt}\bigl( \bkE\bigl| Z_{i-\ell,L}^{K_1,k_{K_1}} \bigr| \bigr) \bigl| \bkE
(Z_{i-j,L}^{K_2,k_{K_2}} \bkE_{i-j+2^L} \bigl(Z_{i,L}^{K_3,k_{K_3}}
\bigr) \bigr|,
\end{eqnarray*}
which implies that
%
\begin{eqnarray}
\label{lind18}\qquad
&&\sum_{j=1}^{u_L-1} \bigl| \bkE
\bigl( D^2 g * \varphi_a (S_{i-1}) -
D^2 g * \varphi_a (S_{i-2j}) \bigr) \blddot\bkE
(Z_{i-j,L} \otimes Z_{i,L} ) \bigr| \I_{j \leq[i/2]}
\nonumber\\[-8pt]\\[-8pt]
&&\qquad\ll a^{-2} \bigl(r(L)\bigr)^2 2^{r(L)}\sum
_{j=1}^{u_L }j \beta_{1,X}(j).
\nonumber
\end{eqnarray}
Therefore, (\ref{lind17}) together with (\ref{lind18}), the choice of
$u_L$ and the condition on the $\beta$-dependence coefficients entail that
%
\begin{equation}
\label{lind19} \qquad\sum_{i=1}^{2^{m(L)}}A_i
\ll a^{-1} L^2 2^{r(L)} + a^{-2}
L^3 2^{r(L)+m(L)} + a^{-1-\delta}L^{1 + \delta}
2^{r(L)+m(L)}.
\end{equation}
Taking into account (\ref{lind1})--(\ref{lind5first}), (\ref{lind8}),
(\ref{lind16}) and (\ref{lind19}), the bound (\ref{boundWasser3})
follows.
\end{pf}

\subsection{\texorpdfstring{Proof of Theorem \protect\ref{maindyn}}{Proof of Theorem 3.1}}
Let $(X_i)_{i \in{\mathbb Z}}$ be a stationary Markov chain with
transition Kernel $Q$
defined in (\ref{Perron}). Notice that
for all $(s,s')\in[0,1]^2$,
\[
\nu\bigl(f_s^{(0)}\cdot f_{s'}^{(0)}
\circ T^k\bigr) = \operatorname{Cov} ( \I_{X_k \leq s},
\I_{X_0 \leq s'}).
\]
Since $\beta_{2,X} (k)$ satisfies (\ref{boundonbeta}), according to
the proof of item (1) of Theorem~\ref{mainemp}, it follows that item (1)
of Theorem~\ref{maindyn} holds true.

As at the beginning of the proof of Theorem~\ref{mainemp}, we start by
considering the probability $P^*_{\nu}$ whose
density with respect to $\nu$ is given by (\ref{defnewprob}). Let
$F^*_{\nu}$ be the distribution function of $P_{\nu}^*$ ($F^*_{\nu}$
is continuous since $\nu$ is absolutely continuous with respect to the
Lebesgue measure). Let now $\widetilde T_i = F^*_{\nu} (T^i)$ and
$Y_i=F^*_{\nu}(X_i)$. Let $F_Y$ be the distribution function of $Y_0$.
Clearly,
$R_{ T}(\cdot, \cdot)= R_{ \widetilde T} (F^*_{\nu}(\cdot), \cdot
)$ almost surely, where
\[
R_{ \widetilde T}(s,t)= \sum_{1 \leq k \leq t} \bigl(
\I_{\widetilde
T_k \leq s}-F_Y(s) \bigr),\qquad s \in{[0, 1]}, t \in{\mathbb
R}^+.
\]
Theorem~\ref{maindyn} will then follow if we can prove that there
exists a two-parameter Gaussian process $K^*_{\widetilde T}$ with
covariance function $\Gamma_{\widetilde T}$ given by $\Gamma
_{\widetilde T} (s,s',t,t') = \min(t,t') \Lambda_{\widetilde T
}(s,s')$, where
%
\begin{equation}
\label{defgammaTtilde}\qquad\quad \Lambda_{\widetilde T} \bigl(s,s'\bigr)=
\sum_{k \geq0} \nu\bigl(f_s^{(0)}
\cdot f_{s'}^{(0)}\circ F^*_{\nu}\bigl(
T^k\bigr)\bigr) + \sum_{k >0} \nu
\bigl(f_{s'}^{(0)}\cdot f_{s}^{(0)}\circ
F^*_{\nu}\bigl( T^k\bigr)\bigr).
\end{equation}

For $L \in{\mathbb N}$, let $m(L)$ and $r(L)$ be the two sequences of
integers defined by (\ref{choicerm}). For any integer $j$, let $s_j =
j 2^{-r(L)}$.
As for the proof of Theorem~\ref{mainemp}, we start by constructing
the approximating Kiefer process $K^*_{\widetilde T}$ with covariance
function $\Gamma_{\widetilde T }$. With this aim, we first define for
any $\ell\in\{1, \ldots, 2^{L-m(L)} \} $,
\[
I_{L, \ell} = {\bigl]2^L + (\ell-1)2^{m(L)},
2^L + \ell2^{m(L)}\bigr]} \cap{\mathbb N}
\]
and
\[
U^{* (j)}_{L,\ell} = \sum_{i\in
I_{L, \ell}}
\bigl(\I_{\widetilde T_i \leq s_j}-F_Y (s_j) \bigr).
\]
The associated column vectors $U^*_{L, \ell}$ are then defined in
${\mathbb R}^{2^{r(L)} -1}$ by the equality
$U^*_{L, \ell} = (U^{* (1)}_{L,\ell}, \ldots, U^{* (2^{r(L)}
-1)}_{L,\ell} )'$.
Let
\[
\Lambda_{\widetilde T,L} = \bigl(\Lambda_{\widetilde T}(s_j,
s_{j'} ) \bigr)_{
j,j'= 1, \ldots, 2^{r(L)} -1 },
\]
where the $ \Lambda_{ \widetilde T} (s_j, s_{j'} )$ are defined in
(\ref{defgammaTtilde}). Let
$G_{2^{m(L)} \Lambda_{\widetilde T,L}}$ denote the ${\mathcal N} (
0,\break
2^{m(L)} \Lambda_{\widetilde T,L})$-law, and for any $\ell\in\{1,
\ldots, 2^{L-m(L)} \} $, let $P_{U^*_{L, \ell} | {\mathcal G}_{2^L
+\ell2^{m(L) }+1
}}$ be the conditional law of $U^*_{L,\ell}$ given ${\mathcal G}_{2^L
+\ell2^{m(L)}+1}$, where $ {\mathcal G}_{m} = \sigma( T^i, {i \geq
m})$. By the Markov property, the
following equality holds: $P_{U^*_{L, \ell} | {\mathcal G}_{2^L +\ell
2^{m(L) }+1
}}=P_{U^*_{L, \ell} | T^{2^L +\ell2^{m(L) }+1}
}$.

According to \citet{Rus85}, there exists $V^*_{L, \ell} = (
V^{*(1)}_{L,\ell}, \ldots, V^{*(2^{r(L)} -1)}_{L,\ell} )'$ with law
$G_{2^{m(L)} \Lambda_{\widetilde T,L}}$, measurable with respect to
$\sigma(\delta_{2^L + \ell2^{m(L)}}) \vee\sigma(U^*_{L,\ell})
\vee{\mathcal G}_{2^L +\ell2^{m(L)} +1}$, independent of ${\mathcal
G}_{2^L +\ell2^{m(L)} +1}$, and such that, with the notation of Section
\ref{constructkiefer},
%
\begin{equation}
\label{coupling1dyn}\quad {\mathbb E} \bigl( d_{r(L)} \bigl(U^*_{L,\ell},
V^*_{L,\ell}\bigr) \bigr) = {\mathbb E} \bigl( W_{d_{r(L)}}
(P_{U^*_{L, \ell} | {\mathcal
G}^*_{2^L + \ell2^{m(L)} +1}},G_{2^{m(L)} \Lambda_{\widetilde T,L}})
\bigr).
\end{equation}
By induction on $\ell$, the random variables $(V^*_{L,\ell})_{\ell= 1,
\ldots, 2^{L-m(L)}}$ are mutually independent, independent of
${\mathcal G}_{2^{L+1}+1}$ and with law ${\mathcal N} ( 0,
2^{m(L)}\Lambda_{\widetilde T,L})$. Hence, we have constructed Gaussian
random variables $(V^*_{L,\ell})_{L \in {\mathbb N}, \ell=1, \ldots,
2^{L-m(L)}}$ that are mutually independent. In addition, according to
Lemma 2.11 of \citet{DudPhi83}, there exists a Kiefer process
$K^*_{\widetilde T}$ with covariance function $\Gamma_{\widetilde T}$
such that for any $L \in {\mathbb N}$, any $\ell\in\{1,\ldots,
2^{L-m(L)} \} $ and any $j \in\{1,\ldots, 2^{r(L) - 1} \}$,
%
\begin{equation}
\label{relationVKieferdyn} V^{*(j)}_{L, \ell} =
K^*_{\widetilde T} \bigl( s_j, 2^L +
\ell2^{m(L)} \bigr) -K^*_{\widetilde T}\bigl(s_j,
2^L + (\ell-1) 2^{m(L)}\bigr).
\end{equation}
Thus, our construction is now complete.

Notice now that, by stationarity, for any $\ell\in\{1, \ldots,
2^{L-m(L)} \} $,
\[
{\mathbb E} \bigl( d_{r(L)} \bigl(U^*_{L,\ell},
V^*_{L,\ell}\bigr) \bigr)={\mathbb E} \bigl( d_{r(L)}
\bigl(U^*_{L,1}, V^*_{L,1}\bigr) \bigr).
\]
In addition, on the probability space $([0, 1], \nu)$, the random
variable $(T^{2^L+1},\break T^{2^L+2},\ldots, T^{2^{L+1}})$ is
distributed as $(X_{2^{L+1}},X_{2^{L+1}-1},\ldots, X_{2^L +1})$. Let
$U^{(j)}_{L,\ell} = \sum_{i\in I_{L, \ell}} (\I_{Y_i \leq
s_j}-F_Y(s_j))$, and
let $U_{L, \ell}$ be the associated column vectors in ${\mathbb
R}^{2^{r(L)} -1}$ defined by $
U_{L, \ell} = (U^{(1)}_{L,\ell}, \ldots, U^{(2^{r(L)}
-1)}_{L,\ell} )'$. According to the coupling relation (\ref
{coupling1}), we get that
%
\begin{eqnarray}
\label{equalitylaw1}
&&{\mathbb E} \bigl( W_{d_{r(L)}} (P_{U^*_{L, 1} |
{\mathcal G}_{2^L +
2^{m(L)} +1}},G_{2^{m(L)} \Lambda_{T,L}})
\bigr)
\nonumber\\
&&\qquad= \bkE\sup_{f
\in\operatorname{Lip} (d_{r(L)} ) } \bigl( {\mathbb E} \bigl( f\bigl(U^*_{L, 1}
\bigr) |T^{2^L + \ell2^{m(L)} +1} \bigr) - {\mathbb E} \bigl(f\bigl
(V^*_{L, 1}
\bigr)\bigr) \bigr)
\\
&&\qquad= \bkE\sup_{f \in\operatorname{Lip} (d_{r(L)} ) } \bigl( {\mathbb E} \bigl(
f(U_{L, 2^{L-m(L)}}) |
X_{2^{L+1}-2^{m(L)}} \bigr) - {\mathbb E} \bigl(f\bigl(V^*_{L, 1}\bigr)
\bigr) \bigr).
\nonumber
\end{eqnarray}
Let us construct the Gaussian random variables $V_{L,\ell}$ associated
to the $U_{L, \ell}$ as in Section~\ref{constructkiefer}. Notice that
since the covariance function $\Lambda_{\widetilde T}$ is the same as
the covariance function $\Lambda_Y$ defined by (\ref{defgammaY}), for
any measurable function $f$, ${\mathbb E} (f(V^*_{L, 1})) ={\mathbb E}
(f(V_{L, 2^{L-m(L)}})) $. Therefore, starting from (\ref
{coupling1dyn}) and taking into account (\ref{equalitylaw1}) together
with (\ref{coupling1}), we get that
%
\begin{eqnarray}
\label{equalitylaw2}\qquad
&&{\mathbb E} \bigl( d_{r(L)} \bigl(U^*_{L,1},
V^*_{L,1}\bigr) \bigr) \nonumber\\
&&\qquad = \bkE\sup_{f \in\operatorname{Lip} (d_{r(L)} ) }
\bigl( {\mathbb
E} \bigl( f(U_{L,
2^{L-m(L)}}) | {\mathcal F}_{2^{L+1}-2^{m(L)}} \bigr) - {\mathbb
E} \bigl(f(V_{L, 2^{L-m(L)}})\bigr) \bigr)
\\
&&\qquad ={\mathbb E} \bigl( d_{r(L)} (U_{L,2^{L-m(L)}}, V_{L,2^{L-m(L)}})
\bigr).\nonumber
\end{eqnarray}
Setting $\Pi_{r(L)}(s) = 2^{-r(L)} [s2^{r(L)}]$ and mimicking the
notation of Section~\ref{sectionerror}, let now
\begin{eqnarray*}
D^*_{L,1} &=& \sup_{2^L < \ell\leq2^{L+1} }\sup_{s \in[0,1]} \bigl|
\bigl(R_{\widetilde T}(s,\ell) -R_{\widetilde T}\bigl(\Pi_{r(L)}(s),
\ell\bigr) \bigr)\\
&&\hspace*{73.1pt}{} - \bigl(R_{\widetilde T}\bigl(s,2^L\bigr)
-R_{\widetilde T}\bigl(\Pi_{r(L)}(s), 2^L \bigr) \bigr) \bigr|,
\\
B^*_{L,3} &=& \sup_{j \in\{ 1,\ldots, 2^{r(L)}-1 \}} \sup_{1 \leq
k \leq2^{L-m(L)}}
\sup_{ \ell\in I_{L,k}} \bigl| R_{\widetilde
T}(s_j, \ell) -
R_{\widetilde T}\bigl(s_j, 2^L + (k-1)2^{m(L)}
\bigr) \bigr|,
\end{eqnarray*}
and let $D_{L,1}$ and $B_{L,3}$ be the same quantities with $R_Y$
replacing $R_{\widetilde T}$. Using once again that, on $([0, 1],
\nu)$, the random variable $(T^{2^L+1}, T^{2^L+2},\ldots,
T^{2^{L+1}})$ is distributed as the random variable
$(X_{2^{L+1}},X_{2^{L+1}-1},\ldots, X_{2^L +1})$, we infer that for
any positive $\lambda$,
%
\begin{equation}
\label{equalitylaw3}\qquad {\mathbb P} \bigl(D^*_{L,1} \geq\lambda\bigr)
\leq{\mathbb P} ( 2 D_{L,1}\geq\lambda)
\quad\mbox{and}\quad {\mathbb P}
\bigl(B^*_{L,3} \geq\lambda\bigr) \leq{\mathbb P} (2 B_{L,3}
\geq\lambda).
\end{equation}
Proceeding as in Section~\ref{sectionerror} of the proof of Theorem
\ref{mainemp}, using the fact that the covariance function $\Gamma
_{\widetilde T}$ is the same as the covariance function $\Gamma_Y$
defined by (\ref{defgammaY}) (so that all the quantities involving
only the Kiefer process $K_{\widetilde T}^*$ can be computed as in
Section~\ref{sectionerror}) and taking into account (\ref
{equalitylaw2}), (\ref{equalitylaw3})
and the fact that the Markov chain $(X_i)_{i \in{\mathbb Z}}$
satisfies the assumptions of Theorem~\ref{mainemp}, Theorem \ref
{maindyn} follows.

\begin{appendix}\label{app}
\section*{Appendix}

\subsection{Properties of the random variables $Y_i$}

For the next lemma, we keep the same notation as that of Definition
\ref{beta} and of the beginning
of Section~\ref{sec4.1}.
Recall that the random variables $Y_i$ have been defined in (\ref{newrv}).
%
\begin{lem}\label{propY} The following assertions hold:
\begin{longlist}[(2)]
\item[(1)] The image measure of ${\mathbb P}_0^*$ by the variable $Y_0$ is
the uniform distribution over $[0,1]$.
%
\item[(2)] The equality
$F_{P^*}^{-1}(Y_i)=X_i$ holds ${\mathbb P}$-almost surely. Moreover,
${\mathbb P}$-almost surely,
\begin{eqnarray*}
b(X_0, k) & \geq& \sup_{t \in{\mathbb R}} \bigl|{P}_{Y_k|X_0}(f_t)-P_Y(f_t)\bigr|,
\\
b_1({\mathcal F}_{\ell}, k) &\geq& \sup_{t \in{\mathbb
R}}\bigl|{P}_{Y_{k}|{\mathcal F}_\ell}(f_t)-P_Y(f_t)\bigr|,
\\
b_2({\mathcal F}_{\ell}, i,j)&\geq&\sup_{(s,t) \in{\mathbb
R}^2}\bigl|P_{(Y_i, Y_j)|{\mathcal F}_\ell}
\bigl(f_t^{(0)}\otimes f_s^{(0)}
\bigr)-P_{(Y_i, Y_j)}\bigl(f_t^{(0)}\otimes
f_s^{(0)}\bigr)\bigr|.
\end{eqnarray*}
\end{longlist}
\end{lem}

\begin{pf}
As in Definition~\ref{beta}, define
\[
b(X_i, k)= \sup_{t \in{\mathbb R}} \bigl|{P}_{X_k|X_i}(f_t)-P(f_t)\bigr|.
\]
On $\Omega$, we introduce the probability ${\mathbb P}^*_i$ whose
density with respect
to ${\mathbb P}$ is
%
\setcounter{equation}{0}
\begin{eqnarray}
\label{defnewProbi} C (\beta)^{-1} \Biggl( 1+4 \sum
_{k=i+1}^\infty b(X_i, k) \Biggr)\nonumber\\[-8pt]\\[-8pt]
&&\eqntext{\displaystyle \mbox{with }C(\beta) = 1+4 \sum_{k=1}^\infty\beta
\bigl(\sigma(X_0), X_k\bigr).}
\end{eqnarray}
By stationarity of $(X_i)_{i \in{\mathbb Z}}$, the image measure of
${\mathbb P}_i^*$ by $X_i$ is again $P^*$.
It follows from Lemma F.1, page 161, in \citet{Rio00} that
the image measure of ${\mathbb P}_i^*$ by the variable $Y_i$ is the
uniform distribution over $[0,1]$ [proving item~(1)],
and that the equality $F_{P^*}^{-1}(Y_i)=X_i$ holds
${\mathbb P}_i^*$-almost surely.
Since the probabilities ${\mathbb P}$ and ${\mathbb P}_i^*$ are equivalent,
it follows that the equality $F_{P^*}^{-1}(Y_i)=X_i$ holds ${\mathbb
P}$-almost surely, proving the first
point of item (2).

Now, note that $Y_i=g(X_i, \eta_i)$, where the function $x \rightarrow
g(x,u)$ is nondecreasing for any $u \in[0,1]$. Since $(X_0, X_k)$ is
independant of $\eta_k$,
\begin{eqnarray*}
&&
\bigl|{P}_{Y_k|X_0}(f_t)-P_Y(f_t)\bigr|\\
&&\qquad= \biggl|
\int_0^1 \bigl\{{\mathbb E}\bigl(
f_t\bigl(g(X_k,u)\bigr)|X_0\bigr)-{
\mathbb E}\bigl( f_t\bigl(g(X_k,u)\bigr)\bigr)\bigr\} \,du
\biggr| \qquad\mbox{almost surely.}
\end{eqnarray*}
The function $x \rightarrow g(x,u)$ being nondecreasing, we infer that
\[
\bigl|{\mathbb E}\bigl( f_t\bigl(g(X_k,u)
\bigr)|X_0\bigr)-{\mathbb E}\bigl( f_t
\bigl(g(X_k,u)\bigr)\bigr)\bigr| \leq b(X_0, k)
\qquad\mbox{almost
surely,}
\]
in such a way that
\[
\bigl|{P}_{Y_k|X_0}(f_t)-P_Y(f_t)\bigr|
\leq b(X_0, k) \qquad\mbox{almost surely.}
\]
The two last inequalities of item (2) may be proved in the same way.
\end{pf}

\subsection{Some upper bounds for partial derivatives}
Let $x$ and $y$ be two column vectors of ${\mathbb R}^{2^{r(L)} -1}$ with
coordinates
\[
x = \bigl( \bigl( x^{(i,k_i)}, k_i\in{\mathcal E}(L,i)
\bigr)_{i =0,\ldots, r(L) -1 } \bigr)'
\]
and
\[
y = \bigl( \bigl(
y^{(i,k_i)}, k_i\in{\mathcal E}(L,i) \bigr)_{i = 0,\ldots, r(L) -1 }
\bigr)',
\]
where ${\mathcal E}(L,i) = \{1, \ldots, 2^{r(L)-i}-1\} \cap(
2{\mathbb N} +1 )$. Let $f \in\operatorname{Lip}(d^*_{r(L)})$, meaning that
\[
\bigl| f(x) - f(y) \bigr| \leq\sum_{K=0}^{r(L)-1}
\sup_{k \in{\mathcal
E}(L,K)} \bigl| x^{(K,k)} -y^{(K,k)} \bigr|
\]
[the distance $d^*_{r(L)}$ is defined in Definition~\ref{defdstar}].
Let $a>0$ and $\varphi_a$ be the density of a centered Gaussian law of
${\mathbb R}^{2^{r(L)} -1}$ with covariance $a^2 \mathrm{I}_{2^{r(L)} -1}
$ ($\mathrm{I}_{2^{r(L)} -1}$ being the identity matrix on ${\mathbb
R}^{2^{r(L)} -1}$). Let also
\[
\| x \|_{\infty,L} = \sum_{K=0}^{r(L) -1}
\sup_{k\in
{\mathcal E}(L,K)} \bigl|x^{(K,k)}\bigr|
\]
and
\[
\| x \|_{2,L}=
\Biggl( \sum_{K=0}^{r(L) -1} \sum
_{k_K\in{\mathcal E}(L,K)} \bigl(x^{(K,k_K)}\bigr)^2
\Biggr)^{1/2}.
\]
For the statements of the lemmas, we refer to Notation~\ref{deftensorprod}.
%
\begin{lem} \label{lma1} The partial derivatives of $f$ exist almost
everywhere and the following inequality holds:
%
\begin{equation}
\label{evidentlip} \sup_{y \in{\mathbb R}^{2^{r(L)} -1}}\sup_{u \in
{\mathbb
R}^{2^{r(L)} -1 }, \| u\|_{\infty,L} \leq1} \bigl| D f (y) \blddot u \bigr|
\leq1.
\end{equation}
In addition,
%
\begin{equation}
\label{inelma1} \sup_{K \in\{ 0, \ldots, r(L) -1 \} } \sum_{k_K\in
{\mathcal
E}(L,K)} \biggl|
\frac{\partial f}{\partial x^{(K,k_K)}} (y) \biggr| \leq1.
\end{equation}
\end{lem}

\begin{pf}
The first part of the lemma
follows directly from the fact that $f$ is Lipschitz with respect to
the distance $d^*_{r(L)}$ together with the Rademacher theorem. We
prove now (\ref{inelma1}). For any $K \in\{0, \ldots, r(L) -1 \}$,
we consider the column vector
$u_K = ( ( u_K^{(i,k_i)}, k_i\in{\mathcal E}(L,i) )_{i
= 0,\ldots, r(L) -1 } )'
$ with coordinates given by
\[
u_K^{(i,k_i)} = \operatorname{sign} \biggl( \frac{\partial f}{\partial
x^{(i,k_i)}} (y)
\biggr) {\mathbf1}_{i=K}.
\]
Applying inequality (\ref{evidentlip}) together with the fact that
$\|u_K \|_{\infty,L}=1$, we get that
\[
\sum_{k\in{\mathcal E}(L,K)} \biggl| \frac{\partial f}{\partial
x^{(K,k)}} (y) \biggr| = \bigl| D f (y)
\blddot u_K \bigr| \leq1
\]
and (\ref{inelma1}) follows.
\end{pf}
%
\begin{lem} \label{lma2} Let $X$ and $Y $ be two random variables in
${\mathbb R}^{2^{r(L)} -1}$. For any positive integer $m $ and any $t
\in[0,1]$,
\[
\bigl| \bkE\bigl( D^m f * \varphi_a (Y +tX) \blddot
X^{\otimes m} \bigr) \bigr| \leq\bkE\bigl( \bigl\| Df(\cdot) \blddot X
\bigr\|_{\infty} \times\bigl\| D^{m-1}\varphi_a (\cdot)
\blddot X^{\otimes
m-1} \bigr\|_1 \bigr).
\]
\end{lem}

\begin{pf}
For any positive integer $m
$ and any $x,y \in{\mathbb R}^{2^{r(L)} -1}$, it follows, from the
properties of the convolution product, that
\[
D^m f * \varphi_a (y) \blddot x^{\otimes m} =
\bigl(Df(\cdot) \blddot x \bigr) * \bigl(D^{m-1}\varphi_a (
\cdot) \blddot x^{\otimes m-1} \bigr) (y),
\]
where $Df(\cdot) \blddot x\dvtx  y \mapsto Df(y) \blddot x$ and
$D^{m-1}\varphi_a (\cdot) \blddot x^{\otimes m-1}\dvtx  y
\mapsto D^{m-1}\varphi_a (y) \blddot x^{\otimes m-1}$. The lemma
then follows immediately.
\end{pf}
%
\begin{lem} \label{lma3} Let $X$ be a random variable in ${\mathbb
R}^{2^{r(L)} -1}$. For any nonnegative integer $m $, there exists a
positive constant $c_m$ depending only on $m$ such that
%
\begin{equation}
\label{majlma3} \bigl\| D^{m}\varphi_a (\cdot) \blddot
X^{\otimes m} \bigr\|_1 \leq c_m a^{-m}
\| X \|_{2,L}^m.
\end{equation}
\end{lem}

\begin{pf}
In order to simplify the
proof, and to avoid the double indexes $(K,k_K)$ for the coordinates of
a column vector of ${\mathbb R}^{2^{r(L)} -1}$, we set $d=2^{r(L)} -1$
and we denote by $x=(x_1,\ldots, x_d)'$ an element of ${\mathbb
R}^{d}$. Proceeding by induction on $m$, we infer that for any $u,x$ in
${\mathbb R}^{d}$ and any integer~$m$,
%
\begin{eqnarray}
\label{formulederivee} D^{m}\varphi_a (u) \blddot
x^{\otimes m} &=& \frac{1}{(2 \pi
a^2)^{d/2}} \nonumber\\[-4pt]\\[-12pt]
&&{}\times\exp\Biggl( -\frac{1}{2a^2} \sum
_{i=1}^d u^2_i
\Biggr) \sum_{\ell=0}^{[m/2]} \frac{c_{m,\ell} }{a^{2 \ell}}
\Biggl( \sum_{i=1}^d x^2_i
\Biggr)^{\ell} \Biggl( \sum_{i=1}^d
\frac{u_i
x_i}{a^2} \Biggr)^{m - 2\ell}\hspace*{-20pt}\nonumber
\end{eqnarray}
with the following recurrence relations between the $c_{m,\ell}$:
\begin{eqnarray*}
c_{m,0}&=&(-1)^m \qquad\mbox{for any }m \geq0,\qquad
c_{2, 1} = -1,
\\
c_{m+1, \ell} &=& ( m- 2 \ell+2) c_{m,\ell-1} - c_{m,\ell}
\end{eqnarray*}
for $\ell\in\bigl\{1,\ldots, [m/2] \bigr\}$ and $m\geq2$,
\begin{eqnarray*}
c_{m+1, [(m+1)/2]} &= & c_{m,[m/2]} \qquad\mbox{if $m $ is odd,}\\
c_{m+1, [(m+1)/2]} &=& c_{m+1, [m/2]} \qquad\mbox{if $m $ is even}.
\end{eqnarray*}
Starting from (\ref{formulederivee}) and setting $\| x \|_{2,d} =
( \sum_{i=1}^d x_i^2 )^{1/2}$, we get that for any integer~$m$,
\begin{eqnarray*}
\hspace*{-4pt}&&\int_{{\mathbb R}^d} \bigl| D^{m}
\varphi_a (u) \blddot x^{\otimes m} \bigr| \,du
\\
\hspace*{-4pt}&&\quad\leq\frac{\| x \|_{2,d}^{m}}{a^m (2 \pi a^2)^{d/2}} \int_{{\mathbb
R}^d} \exp\Biggl( -
\frac{1}{2a^2} \sum_{i=1}^d
u^2_i \Biggr) \sum_{\ell=0}^{m }
\Biggl| c_{m,\ell} \Biggl( \sum_{i=1}^d
\frac{u_i
x_i}{a \| x \|_{2,d} } \Biggr)^{m - 2\ell} \Biggr| \prod_{i=1}^d
d u_i
\\
\hspace*{-4pt}&&\quad\leq\frac{ \| x \|_{2,d}^{m} }{a^{m} } \int_{{\mathbb R}^d}\frac
{1}{(2 \pi)^{d/2} } \exp
\Biggl( -\frac{1}{2} \sum_{i=1}^d
u^2_i \Biggr) \sum_{\ell=0}^{m }
\Biggl| c_{m,\ell} \Biggl( \sum_{i=1}^d
\frac
{u_i x_i}{ \| x \|_{2,d} } \Biggr)^{m - 2\ell} \Biggr| \prod_{i=1}^d
d u_i.
\end{eqnarray*}
Now, for any integer $k$, we have that
\[
\frac{1}{(2 \pi)^{d/2}}\int_{{\mathbb R}^d} \exp\Biggl( -
\frac
{1}{2} \sum_{i=1}^d
u^2_i \Biggr) \Biggl| \sum_{i=1}^d
\frac{u_i x_i}{
\|x\|_{2,d}} \Biggr|^{k} \prod_{i=1}^d
d u_i = \bkE\bigl( |N|^k\bigr),
\]
where $N \sim{\mathcal N} (0,1)$. Therefore,
\[
\int_{{\mathbb R}^d} \bigl| D^{m}\varphi_a (u)
\blddot x^{\otimes m} \bigr| \,du \leq a^{-m}\| x \|_{2,d}^{m}
\sum_{\ell=0}^{[m/2]} |c_{m,\ell}| \bkE
\bigl( |N|^{m-2\ell}\bigr),
\]
which completes the proof of (\ref{majlma3}).
\end{pf}
%
\begin{lem} \label{lma4} Let $X $ and $Y $ be two random variables
with values in ${\mathbb R}^{2^{r(L)} -1}$. For any positive integer $m
$ and any $t \in[0,1]$, there exists a positive constant $c_{m-1}$
depending only on $m$ such that
\[
\bigl| \bkE\bigl( D^m f * \varphi_a (Y +tX) \blddot
X^{\otimes m} \bigr) \bigr| \leq c_{m-1} a^{1-m} \bkE\bigl(
\| X \|_{\infty,L} \times\| X \|^{m-1}_{2,L}
\bigr).
\]
\end{lem}
\begin{pf}
Applying Lemmas~\ref{lma2}
and~\ref{lma3} and using the fact that, by~(\ref{evidentlip}),
\[
\bigl\| Df(\cdot) \blddot X \bigr\|_{\infty} = \| X \|_{\infty,L}
\sup_{y \in{\mathbb R}^{2^{r(L)} -1}} \biggl| Df(y) \blddot\frac{X}{ \| X
\|_{\infty,L} } \biggr| \leq\| X \|_{\infty,L},
\]
the result follows.
\end{pf}
%
\begin{lem} \label{lma5} For any $y \in{\mathbb R}^{2^{r(L)} -1}$ and
any integer $m \geq1$, there exists a positive constant $c_m$
depending only on $m$ such that
\[
\sup_{(K_i, k_{K_i}), i =1,\ldots, m} \biggl| \frac{\partial^m f *
\varphi_a}{\prod_{i=1}^m\partial x^{(K_i,k_{K_i})}} (y) \biggr| \leq c_m
a^{1-m },
\]
where the supremum is taken over all the indexes $K_i \in\{0, \ldots,
r(L) -1 \}$ and $k_{K_i} \in{\mathcal E}(L,K_i)$ for
any $i=1, \ldots, m$.
\end{lem}

\begin{pf}
Notice first that by the
properties of the convolution product,
\[
\frac{\partial^m f * \varphi_a}{\prod_{i=1}^m\partial
x^{(K_i,k_{K_i})}} (y) = \biggl( \frac{\partial f }{\partial
x^{(K_1,k_{K_1})}} * \frac{\partial^{m-1} \varphi_a}{\prod
_{i=2}^m\partial x^{(K_i,k_{K_i})}} \biggr)
(y).
\]
Therefore, by using (\ref{inelma1}),
%
\begin{eqnarray}
\label{p1lma5} \biggl| \frac{\partial^m f * \varphi_a}{\prod_{i=1}^m\partial
x^{(K_i,k_{K_i})}} (y) \biggr| &\leq&\biggl\| \frac{\partial f }{\partial
x^{(K_1,k_{K_1})}}
\biggr\|_{\infty} \biggl\| \frac{\partial^{m-1}
\varphi_a}{\prod_{i=2}^m\partial x^{(K_i,k_{K_i})}} \biggr\|_1
\nonumber\\[-8pt]\\[-8pt]
&\leq&\biggl\|
\frac{\partial^{m-1} \varphi_a}{\prod_{i=2}^m\partial
x^{(K_i,k_{K_i})}} \biggr\|_1.\nonumber
\end{eqnarray}
Let now $h_a$ be the density of the ${\mathcal N} (0,a^2)$
distribution, and let
\[
{\mathcal S}_m = \Biggl\{ (\ell_1,\ldots,
\ell_{m} ) \in\{ 0,\ldots, m \}^{\otimes m} \mbox{ such that }
\sum_{i=1}^m \ell_i = m
\Biggr\}.
\]
With this notation, we infer that
\[
\biggl\| \frac{\partial^{m-1} \varphi_a}{\prod_{i=2}^m\partial
x^{(K_i,k_{K_i})}} \biggr\|_1 \leq\sup_{(\ell_1,\ldots, \ell
_{m-1}) \in{\mathcal S}_{m-1}} \prod
_{i=1}^{m-1} \bigl\| h^{(\ell
_i)}_a
\bigr\|_1,
\]
where $h^{(\ell_i)}_a$ is the $\ell_i$th derivative of $h_a$. Since
for any real $u$, $h^{(\ell_i)}_a (u) = a^{-(\ell_i +1)}h^{(\ell
_i)}_1 (u/a)$, it follows that $\| h^{(\ell_i)}_a \|_1 =
a^{-\ell_i } \| h_1^{(\ell_i)} \|_1$. Therefore,
%
\begin{equation}
\label{p2lma5} \biggl\| \frac{\partial^{m-1} \varphi_a}{\prod_{i=2}^m\partial
x^{(K_i,k_{K_i})}} \biggr\|_1 \leq a^{1-m}
\sup_{(\ell_1,\ldots,
\ell_{m-1}) \in{\mathcal S}_{m-1}} \prod_{i=1}^{m-1} \bigl\|
h_1^{(\ell_i)} \bigr\|_1.
\end{equation}
Starting from (\ref{p1lma5}) and using (\ref{p2lma5}), the lemma is
proved, with
\[
c_m =\sup_{(\ell_1,\ldots, \ell_{m-1}) \in{\mathcal S}_{m-1}} \prod_{i=1}^{m-1}
\bigl\| h_1^{(\ell_i)} \bigr\|_1.
\]
\upqed\end{pf}
\end{appendix}




\printaddresses


\begin{thebibliography}{33}

\bibitem[\protect\citeauthoryear{Berkes, H{\"o}rmann and
  Schauer}{2009}]{BerHorSch09}
\begin{barticle}[mr]
\bauthor{\bsnm{Berkes},~\bfnm{Istv{\'a}n}\binits{I.}},
  \bauthor{\bsnm{H{\"o}rmann},~\bfnm{Siegfried}\binits{S.}} \AND
  \bauthor{\bsnm{Schauer},~\bfnm{Johannes}\binits{J.}}
(\byear{2009}).
\btitle{Asymptotic results for the empirical process of stationary sequences}.
\bjournal{Stochastic Process. Appl.}
\bvolume{119}
\bpages{1298--1324}.
\bid{doi={10.1016/j.spa.2008.06.010}, issn={0304-4149}, mr={2508575}}
\bptok{imsref}%
\end{barticle}
\endbibitem

\bibitem[\protect\citeauthoryear{Berkes and Philipp}{1977}]{BerPhi77}
\begin{barticle}[mr]
\bauthor{\bsnm{Berkes},~\bfnm{Istv{\'a}n}\binits{I.}} \AND
  \bauthor{\bsnm{Philipp},~\bfnm{Walter}\binits{W.}}
(\byear{1977}).
\btitle{An almost sure invariance principle for the empirical distribution
  function of mixing random variables}.
\bjournal{Z. Wahrsch. Verw. Gebiete}
\bvolume{41}
\bpages{115--137}.
\bid{mr={0464344}}
\bptok{imsref}%
\end{barticle}
\endbibitem

\bibitem[\protect\citeauthoryear{Bickel and Wichura}{1971}]{BicWic71}
\begin{barticle}[mr]
\bauthor{\bsnm{Bickel},~\bfnm{P.~J.}\binits{P.~J.}} \AND
  \bauthor{\bsnm{Wichura},~\bfnm{M.~J.}\binits{M.~J.}}
(\byear{1971}).
\btitle{Convergence criteria for multiparameter stochastic processes and some
  applications}.
\bjournal{Ann. Math. Statist.}
\bvolume{42}
\bpages{1656--1670}.
\bid{issn={0003-4851}, mr={0383482}}
\bptok{imsref}%
\end{barticle}
\endbibitem

\bibitem[\protect\citeauthoryear{Borovkova, Burton and
  Dehling}{2001}]{BorBurDeh01}
\begin{barticle}[mr]
\bauthor{\bsnm{Borovkova},~\bfnm{Svetlana}\binits{S.}},
  \bauthor{\bsnm{Burton},~\bfnm{Robert}\binits{R.}} \AND
  \bauthor{\bsnm{Dehling},~\bfnm{Herold}\binits{H.}}
(\byear{2001}).
\btitle{Limit theorems for functionals of mixing processes with applications to
  {$U$}-statistics and dimension estimation}.
\bjournal{Trans. Amer. Math. Soc.}
\bvolume{353}
\bpages{4261--4318}.
\bid{doi={10.1090/S0002-9947-01-02819-7}, issn={0002-9947}, mr={1851171}}
\bptok{imsref}%
\end{barticle}
\endbibitem

\bibitem[\protect\citeauthoryear{Castelle and
  Laurent-Bonvalot}{1998}]{CasLau98}
\begin{barticle}[mr]
\bauthor{\bsnm{Castelle},~\bfnm{Nathalie}\binits{N.}} \AND
  \bauthor{\bsnm{Laurent-Bonvalot},~\bfnm{Fran{\c{c}}oise}\binits{F.}}
(\byear{1998}).
\btitle{Strong approximations of bivariate uniform empirical processes}.
\bjournal{Ann. Inst. Henri Poincar\'e Probab. Stat.}
\bvolume{34}
\bpages{425--480}.
\bid{doi={10.1016/S0246-0203(98)80024-1}, issn={0246-0203}, mr={1632841}}
\bptok{imsref}%
\end{barticle}
\endbibitem

\bibitem[\protect\citeauthoryear{Dedecker}{2010}]{Ded10}
\begin{barticle}[mr]
\bauthor{\bsnm{Dedecker},~\bfnm{J.}\binits{J.}}
(\byear{2010}).
\btitle{An empirical central limit theorem for intermittent maps}.
\bjournal{Probab. Theory Related Fields}
\bvolume{148}
\bpages{177--195}.
\bid{doi={10.1007/s00440-009-0227-5}, issn={0178-8051}, mr={2653226}}
\bptok{imsref}%
\end{barticle}
\endbibitem

\bibitem[\protect\citeauthoryear{Dedecker, Gou{\"e}zel and
  Merlev{\`e}de}{2010}]{DedGouMer10}
\begin{barticle}[mr]
\bauthor{\bsnm{Dedecker},~\bfnm{J.}\binits{J.}},
  \bauthor{\bsnm{Gou{\"e}zel},~\bfnm{S.}\binits{S.}} \AND
  \bauthor{\bsnm{Merlev{\`e}de},~\bfnm{F.}\binits{F.}}
(\byear{2010}).
\btitle{Some almost sure results for unbounded functions of intermittent maps
  and their associated {M}arkov chains}.
\bjournal{Ann. Inst. Henri Poincar\'e Probab. Stat.}
\bvolume{46}
\bpages{796--821}.
\bid{doi={10.1214/09-AIHP343}, issn={0246-0203}, mr={2682267}}
\bptok{imsref}%
\end{barticle}
\endbibitem

\bibitem[\protect\citeauthoryear{Dedecker and Merlev{\`e}de}{2010}]{DedMer10}
\begin{bincollection}[mr]
\bauthor{\bsnm{Dedecker},~\bfnm{J{\'e}r{\^o}me}\binits{J.}} \AND
  \bauthor{\bsnm{Merlev{\`e}de},~\bfnm{Florence}\binits{F.}}
(\byear{2010}).
\btitle{On the almost sure invariance principle for stationary sequences of
  {H}ilbert-valued random variables}.
In \bbooktitle{Dependence in Probability, Analysis and Number Theory}
\bpages{157--175}.
\bpublisher{Kendrick Press}, \baddress{Heber City, UT}.
\bid{mr={2731073}}
\bptok{imsref}%
\end{bincollection}
\endbibitem

\bibitem[\protect\citeauthoryear{Dedecker, Prieur and Raynaud
  De~Fitte}{2006}]{DedPriRay06}
\begin{bincollection}[mr]
\bauthor{\bsnm{Dedecker},~\bfnm{J{\'e}r{\^o}me}\binits{J.}},
  \bauthor{\bsnm{Prieur},~\bfnm{Cl{\'e}mentine}\binits{C.}} \AND
  \bauthor{\bsnm{Raynaud De~Fitte},~\bfnm{Paul}\binits{P.}}
(\byear{2006}).
\btitle{Parametrized {K}antorovich-{R}ubin\v stein theorem and application to
  the coupling of random variables}.
In \bbooktitle{Dependence in Probability and Statistics}.
\bseries{Lecture Notes in Statistics}
\bvolume{187}
\bpages{105--121}.
\bpublisher{Springer}, \blocation{New York}.
\bid{doi={10.1007/0-387-36062-X_5}, mr={2283252}}
\bptok{imsref}%
\end{bincollection}
\endbibitem

\bibitem[\protect\citeauthoryear{Dedecker and Prieur}{2007}]{DedPri07}
\begin{barticle}[mr]
\bauthor{\bsnm{Dedecker},~\bfnm{J{\'e}r{\^o}me}\binits{J.}} \AND
  \bauthor{\bsnm{Prieur},~\bfnm{Cl{\'e}mentine}\binits{C.}}
(\byear{2007}).
\btitle{An empirical central limit theorem for dependent sequences}.
\bjournal{Stochastic Process. Appl.}
\bvolume{117}
\bpages{121--142}.
\bid{doi={10.1016/j.spa.2006.06.003}, issn={0304-4149}, mr={2287106}}
\bptok{imsref}%
\end{barticle}
\endbibitem

\bibitem[\protect\citeauthoryear{Dedecker and Prieur}{2009}]{DedPri09}
\begin{barticle}[mr]
\bauthor{\bsnm{Dedecker},~\bfnm{J.}\binits{J.}} \AND
  \bauthor{\bsnm{Prieur},~\bfnm{C.}\binits{C.}}
(\byear{2009}).
\btitle{Some unbounded functions of intermittent maps for which the central
  limit theorem holds}.
\bjournal{ALEA Lat. Am. J. Probab. Math. Stat.}
\bvolume{5}
\bpages{29--45}.
\bid{issn={1980-0436}, mr={2475605}}
\bptok{imsref}%
\end{barticle}
\endbibitem

\bibitem[\protect\citeauthoryear{Dedecker et~al.}{2007}]{Dedetal07}
\begin{bbook}[mr]
\bauthor{\bsnm{Dedecker},~\bfnm{J{\'e}r{\^o}me}\binits{J.}},
  \bauthor{\bsnm{Doukhan},~\bfnm{Paul}\binits{P.}},
  \bauthor{\bsnm{Lang},~\bfnm{Gabriel}\binits{G.}},
  \bauthor{\bsnm{Le{\'o}n~R.},~\bfnm{Jos{\'e}~Rafael}\binits{J.~R.}},
  \bauthor{\bsnm{Louhichi},~\bfnm{Sana}\binits{S.}} \AND
  \bauthor{\bsnm{Prieur},~\bfnm{Cl{\'e}mentine}\binits{C.}}
(\byear{2007}).
\btitle{Weak Dependence: With Examples and Applications}.
\bseries{Lecture Notes in Statistics}
\bvolume{190}.
\bpublisher{Springer}, \blocation{New York}.
\bid{mr={2338725}}
\bptok{imsref}%
\end{bbook}
\endbibitem

\bibitem[\protect\citeauthoryear{Dehling and Taqqu}{1989}]{DehTaq89}
\begin{barticle}[mr]
\bauthor{\bsnm{Dehling},~\bfnm{Herold}\binits{H.}} \AND
  \bauthor{\bsnm{Taqqu},~\bfnm{Murad~S.}\binits{M.~S.}}
(\byear{1989}).
\btitle{The empirical process of some long-range dependent sequences with an
  application to {$U$}-statistics}.
\bjournal{Ann. Statist.}
\bvolume{17}
\bpages{1767--1783}.
\bid{doi={10.1214/aos/1176347394}, issn={0090-5364}, mr={1026312}}
\bptok{imsref}%
\end{barticle}
\endbibitem

\bibitem[\protect\citeauthoryear{Dudley and Philipp}{1983}]{DudPhi83}
\begin{barticle}[mr]
\bauthor{\bsnm{Dudley},~\bfnm{R.~M.}\binits{R.~M.}} \AND
  \bauthor{\bsnm{Philipp},~\bfnm{Walter}\binits{W.}}
(\byear{1983}).
\btitle{Invariance principles for sums of {B}anach space valued random elements
  and empirical processes}.
\bjournal{Z. Wahrsch. Verw. Gebiete}
\bvolume{62}
\bpages{509--552}.
\bid{doi={10.1007/BF00534202}, issn={0044-3719}, mr={0690575}}
\bptok{imsref}%
\end{barticle}
\endbibitem

\bibitem[\protect\citeauthoryear{Finkelstein}{1971}]{Fin71}
\begin{barticle}[mr]
\bauthor{\bsnm{Finkelstein},~\bfnm{Helen}\binits{H.}}
(\byear{1971}).
\btitle{The law of the iterated logarithm for empirical distributions}.
\bjournal{Ann. Math. Statist.}
\bvolume{42}
\bpages{607--615}.
\bid{issn={0003-4851}, mr={0287600}}
\bptok{imsref}%
\end{barticle}
\endbibitem

\bibitem[\protect\citeauthoryear{Giraitis and Surgailis}{2002}]{GirSur02}
\begin{bincollection}[mr]
\bauthor{\bsnm{Giraitis},~\bfnm{Liudas}\binits{L.}} \AND
  \bauthor{\bsnm{Surgailis},~\bfnm{Donatas}\binits{D.}}
(\byear{2002}).
\btitle{The reduction principle for the empirical process of a long memory
  linear process}.
In \bbooktitle{Empirical Process Techniques for Dependent Data}
\bpages{241--255}.
\bpublisher{Birkh\"auser}, \blocation{Boston, MA}.
\bid{mr={1958784}}
\bptok{imsref}%
\end{bincollection}
\endbibitem

\bibitem[\protect\citeauthoryear{Hennion and Herv{\'e}}{2001}]{HenHer01}
\begin{bbook}[mr]
\bauthor{\bsnm{Hennion},~\bfnm{Hubert}\binits{H.}} \AND
  \bauthor{\bsnm{Herv{\'e}},~\bfnm{Lo{\"{\i}}c}\binits{L.}}
(\byear{2001}).
\btitle{Limit Theorems for {M}arkov Chains and Stochastic Properties of
  Dynamical Systems by Quasi-Compactness}.
\bseries{Lecture Notes in Math.}
\bvolume{1766}.
\bpublisher{Springer}, \blocation{Berlin}.
\bid{doi={10.1007/b87874}, mr={1862393}}
\bptok{imsref}%
\end{bbook}
\endbibitem

\bibitem[\protect\citeauthoryear{Kiefer}{1972}]{Kie72}
\begin{barticle}[auto:STB|2013/01/29|08:09:18]
\bauthor{\bsnm{Kiefer},~\bfnm{J.}\binits{J.}}
(\byear{1972}).
\btitle{Skorohod embedding of multivariate rv's, and the sample df}.
\bjournal{Z.~Wahrsch. Verw. Gebiete}
\bvolume{24}
\bpages{1--35}.
\bid{mr={0341636}}
\bptok{imsref}%
\end{barticle}
\endbibitem

\bibitem[\protect\citeauthoryear{Koml{\'o}s, Major and
  Tusn{\'a}dy}{1975}]{KomMajTus75}
\begin{barticle}[mr]
\bauthor{\bsnm{Koml{\'o}s},~\bfnm{J.}\binits{J.}},
  \bauthor{\bsnm{Major},~\bfnm{P.}\binits{P.}} \AND
  \bauthor{\bsnm{Tusn{\'a}dy},~\bfnm{G.}\binits{G.}}
(\byear{1975}).
\btitle{An approximation of partial sums of independent {${\rm RV}$}'s and the
  sample {${\rm DF}$}. {I}}.
\bjournal{Z. Wahrsch. Verw. Gebiete}
\bvolume{32}
\bpages{111--131}.
\bid{mr={0375412}}
\bptok{imsref}%
\end{barticle}
\endbibitem

\bibitem[\protect\citeauthoryear{Lai}{1974}]{Lai74}
\begin{barticle}[mr]
\bauthor{\bsnm{Lai},~\bfnm{Tze~Leung}\binits{T.~L.}}
(\byear{1974}).
\btitle{Reproducing kernel {H}ilbert spaces and the law of the iterated
  logarithm for {G}aussian processes}.
\bjournal{Z. Wahrsch. Verw. Gebiete}
\bvolume{29}
\bpages{7--19}.
\bid{mr={0368121}}
\bptok{imsref}%
\end{barticle}
\endbibitem

\bibitem[\protect\citeauthoryear{Ledoux and Talagrand}{1991}]{LedTal91}
\begin{bbook}[mr]
\bauthor{\bsnm{Ledoux},~\bfnm{Michel}\binits{M.}} \AND
  \bauthor{\bsnm{Talagrand},~\bfnm{Michel}\binits{M.}}
(\byear{1991}).
\btitle{Probability in {B}anach Spaces: Isoperimetry and Processes}.
\bseries{Ergebnisse der Mathematik und Ihrer Grenzgebiete (3) [Results in
  Mathematics and Related Areas (3)]}
\bvolume{23}.
\bpublisher{Springer}, \blocation{Berlin}.
\bid{mr={1102015}}
\bptok{imsref}%
\end{bbook}
\endbibitem

\bibitem[\protect\citeauthoryear{Liverani, Saussol and
  Vaienti}{1999}]{LivSauVai99}
\begin{barticle}[mr]
\bauthor{\bsnm{Liverani},~\bfnm{Carlangelo}\binits{C.}},
  \bauthor{\bsnm{Saussol},~\bfnm{Beno{\^{\i}}t}\binits{B.}} \AND
  \bauthor{\bsnm{Vaienti},~\bfnm{Sandro}\binits{S.}}
(\byear{1999}).
\btitle{A probabilistic approach to intermittency}.
\bjournal{Ergodic Theory Dynam. Systems}
\bvolume{19}
\bpages{671--685}.
\bid{doi={10.1017/S0143385799133856}, issn={0143-3857}, mr={1695915}}
\bptok{imsref}%
\end{barticle}
\endbibitem

\bibitem[\protect\citeauthoryear{Merlev{\`e}de and Rio}{2012}]{MerRio12}
\begin{barticle}[mr]
\bauthor{\bsnm{Merlev{\`e}de},~\bfnm{Florence}\binits{F.}} \AND
  \bauthor{\bsnm{Rio},~\bfnm{Emmanuel}\binits{E.}}
(\byear{2012}).
\btitle{Strong approximation of partial sums under dependence conditions with
  application to dynamical systems}.
\bjournal{Stochastic Process. Appl.}
\bvolume{122}
\bpages{386--417}.
\bid{doi={10.1016/j.spa.2011.08.012}, issn={0304-4149}, mr={2860454}}
\bptok{imsref}%
\end{barticle}
\endbibitem

\bibitem[\protect\citeauthoryear{Rio}{2000}]{Rio00}
\begin{bbook}[mr]
\bauthor{\bsnm{Rio},~\bfnm{Emmanuel}\binits{E.}}
(\byear{2000}).
\btitle{Th\'eorie Asymptotique des Processus Al\'eatoires Faiblement
  D\'ependants}.
\bseries{Math\'ematiques \& Applications (Berlin) [Mathematics \&
  Applications]}
\bvolume{31}.
\bpublisher{Springer}, \blocation{Berlin}.
\bid{mr={2117923}}
\bptok{imsref}%
\end{bbook}
\endbibitem

\bibitem[\protect\citeauthoryear{Rosenblatt}{1956}]{Ros56}
\begin{barticle}[mr]
\bauthor{\bsnm{Rosenblatt},~\bfnm{M.}\binits{M.}}
(\byear{1956}).
\btitle{A central limit theorem and a strong mixing condition}.
\bjournal{Proc. Natl. Acad. Sci. USA}
\bvolume{42}
\bpages{43--47}.
\bid{issn={0027-8424}, mr={0074711}}
\bptok{imsref}%
\end{barticle}
\endbibitem

\bibitem[\protect\citeauthoryear{R{\"u}schendorf}{1985}]{Rus85}
\begin{barticle}[mr]
\bauthor{\bsnm{R{\"u}schendorf},~\bfnm{Ludger}\binits{L.}}
(\byear{1985}).
\btitle{The {W}asserstein distance and approximation theorems}.
\bjournal{Z.~Wahrsch. Verw. Gebiete}
\bvolume{70}
\bpages{117--129}.
\bid{doi={10.1007/BF00532240}, issn={0044-3719}, mr={0795791}}
\bptok{imsref}%
\end{barticle}
\endbibitem

\bibitem[\protect\citeauthoryear{Shorack and Wellner}{1986}]{ShoWel86}
\begin{bbook}[mr]
\bauthor{\bsnm{Shorack},~\bfnm{Galen~R.}\binits{G.~R.}} \AND
  \bauthor{\bsnm{Wellner},~\bfnm{Jon~A.}\binits{J.~A.}}
(\byear{1986}).
\btitle{Empirical Processes with Applications to Statistics}.
\bpublisher{Wiley}, \blocation{New York}.
\bid{mr={0838963}}
\bptok{imsref}%
\end{bbook}
\endbibitem

\bibitem[\protect\citeauthoryear{Wu}{2007}]{Wu07}
\begin{barticle}[mr]
\bauthor{\bsnm{Wu},~\bfnm{Wei~Biao}\binits{W.~B.}}
(\byear{2007}).
\btitle{Strong invariance principles for dependent random variables}.
\bjournal{Ann. Probab.}
\bvolume{35}
\bpages{2294--2320}.
\bid{doi={10.1214/009117907000000060}, issn={0091-1798}, mr={2353389}}
\bptok{imsref}%
\end{barticle}
\endbibitem

\bibitem[\protect\citeauthoryear{Wu}{2008}]{Wu08}
\begin{barticle}[mr]
\bauthor{\bsnm{Wu},~\bfnm{Wei~Biao}\binits{W.~B.}}
(\byear{2008}).
\btitle{Empirical processes of stationary sequences}.
\bjournal{Statist. Sinica}
\bvolume{18}
\bpages{313--333}.
\bid{issn={1017-0405}, mr={2384990}}
\bptok{imsref}%
\end{barticle}
\endbibitem

\bibitem[\protect\citeauthoryear{Yoshihara}{1979}]{Yos79}
\begin{barticle}[mr]
\bauthor{\bsnm{Yoshihara},~\bfnm{Ken-ichi}\binits{K.-i.}}
(\byear{1979}).
\btitle{Note on an almost sure invariance principle for some empirical
  processes}.
\bjournal{Yokohama Math. J.}
\bvolume{27}
\bpages{105--110}.
\bid{issn={0044-0523}, mr={0560618}}
\bptok{imsref}%
\end{barticle}
\endbibitem

\bibitem[\protect\citeauthoryear{Yu}{1993}]{Yu93}
\begin{barticle}[mr]
\bauthor{\bsnm{Yu},~\bfnm{Hao}\binits{H.}}
(\byear{1993}).
\btitle{A {G}livenko--{C}antelli lemma and weak convergence for empirical
  processes of associated sequences}.
\bjournal{Probab. Theory Related Fields}
\bvolume{95}
\bpages{357--370}.
\bid{doi={10.1007/BF01192169}, issn={0178-8051}, mr={1213196}}
\bptok{imsref}%
\end{barticle}
\endbibitem

\bibitem[\protect\citeauthoryear{Zweim{\"u}ller}{1998}]{Zwe98}
\begin{barticle}[mr]
\bauthor{\bsnm{Zweim{\"u}ller},~\bfnm{Roland}\binits{R.}}
(\byear{1998}).
\btitle{Ergodic structure and invariant densities of non-{M}arkovian interval
  maps with indifferent fixed points}.
\bjournal{Nonlinearity}
\bvolume{11}
\bpages{1263--1276}.
\bid{doi={10.1088/0951-7715/11/5/005}, issn={0951-7715}, mr={1644385}}
\bptok{imsref}%
\end{barticle}
\endbibitem

\end{thebibliography}
\end{document}